\documentclass{gtart}

\def\ifplaintex{\expandafter\ifx\csname documentclass\endcsname\relax}

\def\gtp{{\mathsurround=0pt\it $\cal G\mskip-2mu$eometry \&\ 
$\cal T\!\!$opology $\cal P\!$ublications}}  

\def\recd{{\small Received:\qua\receiveddate\ifx\reviseddate\relax
\else\qquad Revised:\qua\reviseddate\fi\par}} 

\def\lognumber#1{\def\thelognumber{#1}}
\def\volumenumber#1{\def\thevolumenumber{#1}}
\def\volumeyear#1{\def\thevolumeyear{#1}}
\def\papernumber#1{\def\thepapernumber{#1}}
\def\pagenumbers#1#2{\def\startpage{#1}\def\finishpage{#2}}
\def\published#1{\def\publishdate{#1}}

\def\received#1{\def\receiveddate{#1}}

\def\accepted#1{\def\accepteddate{#1}}

\def\asciiaddress#1{\def\theasciiaddress{#1}}

\long\def\asciiabstract#1{\long\def\theasciiabstract{#1}}

\let\\\par\let\thelognumber\relax\let\thevolumenumber\relax
\let\thepapernumber\relax\let\thevolumeyear\relax\let\startpage\relax
\let\finishpage\relax\let\publishdate\relax\let\receiveddate\relax
\let\reviseddate\relax\let\accepteddate\relax\let\theasciititle\relax
\let\theasciiauthors\relax\let\theasciiaddress\relax
\let\theasciiabstract\relax

\let\theasciiemail\relax

\ifplaintex
\font\logobig=cmssbx10 scaled 3836
\font\logomed=cmssbx10 scaled 2557
\else
\font\logobig=cmssbx10 scaled 4200
\font\logomed=cmssbx10 scaled 2800
\fi

\long\def\makeagttitle{   
\count0=\startpage
\agt\hfill      
\hbox to 45truept{\vbox to 0pt{\vglue -13truept{\logomed A\kern -.37em{\logobig 
T}\kern -.38em G}\vss}\hss}
\break
{\small Volume \thevolumenumber\ (\thevolumeyear)
\startpage--\finishpage\nl
Published: \publishdate}

\vglue .25truein

{\parskip=0pt\leftskip 0pt plus
1fil\def\\{\par\smallskip}{\Large\bf\thetitle}\par\medskip} \vglue
0.05truein

{\parskip=0pt\leftskip 0pt plus 1fil\def\\{\par}{\sc\theauthors}
\par\medskip}%
 
\vglue 0.03truein 

{\small\leftskip 25truept\rightskip 25truept{\bf Abstract}\stdspace\theabstract

{\bf AMS Classification}\stdspace\theprimaryclass
\ifx\thesecondaryclass\relax\else; \thesecondaryclass\fi\par
{\bf Keywords}\stdspace \thekeywords\par}\vglue 7truept

}   

\ifplaintex
\font\phead=cmsl9 scaled 950
\font\pnum=cmbx10 scaled 913
\font\pfoot=cmsl9 scaled 950
\headline{\vbox to 0pt{\vskip -4.5mm\line{\small\phead\ifnum
\count0=\startpage ISSN 1472-2739 (on-line) 1472-2747 (printed)
\hfill {\pnum\folio}\else\ifodd\count0\def\\{ }%
\ifx\theshorttitle\relax\thetitle\else\theshorttitle\fi\hfill{\pnum\folio}
\else\def\\{ and }{\pnum\folio}\hfill\ifx\theshortauthors\relax\theauthors
\else\theshortauthors\fi\fi\fi}\vss}}
\footline{\vbox to 0pt{\vglue 0mm\line{\small\pfoot\ifnum\count0=\startpage
\copyright\ \gtp\hfill\else
\agt, Volume \thevolumenumber\ (\thevolumeyear)\hfill\fi}\vss}}
\else
\font=cmsl9 scaled 1050
\font\lnum=cmbx10 
\font=cmsl9 scaled 1050
\makeatletter
\def\@oddhead{{\small\lhead\ifnum\count0=\startpage ISSN 1472-2739 
(on-line) 1472-2747 (printed)\hfill {\lnum\number\count0}\else\ifodd\count0
\def\\{ }\ifx\theshorttitle\relax \thetitle \else\theshorttitle\fi\hfill
{\lnum\number\count0}\else\def\\{ and }{\lnum\number\count0}
\hfill\ifx\theshortauthors\relax 
\theauthors\else\theshortauthors\fi\fi\fi}}\def\@evenhead{\@oddhead}
\def\@oddfoot{\small\lfoot\ifnum\count0=\startpage\copyright\ \gtp\hfill\else
\agt, Volume \thevolumenumber\ (\thevolumeyear)\hfill\fi}
\def\@evenfoot{\@oddfoot}
\makeatother
\fi
\let\maketitlepage\makeagttitle
\let\makeshorttitle\maketitlepage
\let\maketitle\maketitlepage

\newwrite\gtoutfile
\long\gdef\makeheadfile{  
{\def\\{, }\def\s{ }
\immediate\openout\gtoutfile head.xxx
\immediate\write\gtoutfile{To: math@arxiv.org}
\immediate\write\gtoutfile{Subject: put OR rep NNNNN:ppppp}
\immediate\write\gtoutfile{--text follows this line--}
\immediate\write\gtoutfile{Proxy-for: \ifx\theasciiauthors\relax
\theauthors\else\theasciiauthors\fi\s<\ifx\theasciiemail\relax\theemail\else\theasciiemail\fi>}
\immediate\write\gtoutfile{\noexpand\\}
\immediate\write\gtoutfile{Authors: \ifx\theasciiauthors\relax
\theauthors\else\theasciiauthors\fi}
{\def\\{ }\immediate\write\gtoutfile{Title: \ifx\theasciititle\relax
\thetitle\else\theasciititle\fi}}
\immediate\write\gtoutfile{Subj-class: GT or SG, GR etc}
\immediate\write\gtoutfile{MSC-class: \theprimaryclass\ifx\thesecondaryclass\relax\else, \thesecondaryclass\fi}
\immediate\write\gtoutfile{Journal-ref: Algebr. Geom. Topol. \thevolumenumber\s
(\thevolumeyear) \startpage-\finishpage}
\immediate\write\gtoutfile{Comments: Published by Algebraic and
Geometric Topology at}
\immediate\write\gtoutfile{\s\s\s  http://www.maths.warwick.ac.uk/agt/AGTVol\thevolumenumber/agt-\thevolumenumber-\thepapernumber.abs.html}
\immediate\write\gtoutfile{\noexpand\\}
\immediate\write\gtoutfile{}
\ifx\theasciiabstract\relax
\immediate\write\gtoutfile{\theabstract}\else
\immediate\write\gtoutfile{\theasciiabstract}\fi
\immediate\write\gtoutfile{}
\immediate\write\gtoutfile{\noexpand\\}
\immediate\write\gtoutfile{}
\immediate\closeout\gtoutfile}}  

\def\maketitlepage{\makeagttitle\makeheadfile}
\let\makeshorttitle\maketitlepage
\let\maketitle\maketitlepage

\def\ifplaintex{\expandafter\ifx\csname documentclass\endcsname\relax}

\def\gtp{{\mathsurround=0pt\it $\cal G\mskip-2mu$eometry \&\ 
$\cal T\!\!$opology $\cal P\!$ublications}}  

\def\recd{{\small Received:\qua\receiveddate\ifx\reviseddate\relax
\else\qquad Revised:\qua\reviseddate\fi\par}} 

\def\lognumber#1{\def\thelognumber{#1}}
\def\volumenumber#1{\def\thevolumenumber{#1}}
\def\volumeyear#1{\def\thevolumeyear{#1}}
\def\papernumber#1{\def\thepapernumber{#1}}
\def\pagenumbers#1#2{\def\startpage{#1}\def\finishpage{#2}}
\def\published#1{\def\publishdate{#1}}

\def\received#1{\def\receiveddate{#1}}

\def\accepted#1{\def\accepteddate{#1}}

\def\asciiaddress#1{\def\theasciiaddress{#1}}

\long\def\asciiabstract#1{\long\def\theasciiabstract{#1}}

\let\\\par\let\thelognumber\relax\let\thevolumenumber\relax
\let\thepapernumber\relax\let\thevolumeyear\relax\let\startpage\relax
\let\finishpage\relax\let\publishdate\relax\let\receiveddate\relax
\let\reviseddate\relax\let\accepteddate\relax\let\theasciititle\relax
\let\theasciiauthors\relax\let\theasciiaddress\relax
\let\theasciiabstract\relax

\let\theasciiemail\relax

\ifplaintex
\font\logobig=cmssbx10 scaled 3836
\font\logomed=cmssbx10 scaled 2557
\else
\font\logobig=cmssbx10 scaled 4200
\font\logomed=cmssbx10 scaled 2800
\fi

\long\def\makeagttitle{   
\count0=\startpage
\agt\hfill      
\hbox to 45truept{\vbox to 0pt{\vglue -13truept{\logomed A\kern -.37em{\logobig 
T}\kern -.38em G}\vss}\hss}
\break
{\small Volume \thevolumenumber\ (\thevolumeyear)
\startpage--\finishpage\nl
Published: \publishdate}

\vglue .25truein

{\parskip=0pt\leftskip 0pt plus
1fil\def\\{\par\smallskip}{\Large\bf\thetitle}\par\medskip} \vglue
0.05truein

{\parskip=0pt\leftskip 0pt plus 1fil\def\\{\par}{\sc\theauthors}
\par\medskip}%
 
\vglue 0.03truein 

{\small\leftskip 25truept\rightskip 25truept{\bf Abstract}\stdspace\theabstract

{\bf AMS Classification}\stdspace\theprimaryclass
\ifx\thesecondaryclass\relax\else; \thesecondaryclass\fi\par
{\bf Keywords}\stdspace \thekeywords\par}\vglue 7truept

}   

\ifplaintex
\font\phead=cmsl9 scaled 950
\font\pnum=cmbx10 scaled 913
\font\pfoot=cmsl9 scaled 950
\headline{\vbox to 0pt{\vskip -4.5mm\line{\small\phead\ifnum
\count0=\startpage ISSN 1472-2739 (on-line) 1472-2747 (printed)
\hfill {\pnum\folio}\else\ifodd\count0\def\\{ }%
\ifx\theshorttitle\relax\thetitle\else\theshorttitle\fi\hfill{\pnum\folio}
\else\def\\{ and }{\pnum\folio}\hfill\ifx\theshortauthors\relax\theauthors
\else\theshortauthors\fi\fi\fi}\vss}}
\footline{\vbox to 0pt{\vglue 0mm\line{\small\pfoot\ifnum\count0=\startpage
\copyright\ \gtp\hfill\else
\agt, Volume \thevolumenumber\ (\thevolumeyear)\hfill\fi}\vss}}
\else
\font=cmsl9 scaled 1050
\font\lnum=cmbx10 
\font=cmsl9 scaled 1050
\makeatletter
\def\@oddhead{{\small\lhead\ifnum\count0=\startpage ISSN 1472-2739 
(on-line) 1472-2747 (printed)\hfill {\lnum\number\count0}\else\ifodd\count0
\def\\{ }\ifx\theshorttitle\relax \thetitle \else\theshorttitle\fi\hfill
{\lnum\number\count0}\else\def\\{ and }{\lnum\number\count0}
\hfill\ifx\theshortauthors\relax 
\theauthors\else\theshortauthors\fi\fi\fi}}\def\@evenhead{\@oddhead}
\def\@oddfoot{\small\lfoot\ifnum\count0=\startpage\copyright\ \gtp\hfill\else
\agt, Volume \thevolumenumber\ (\thevolumeyear)\hfill\fi}
\def\@evenfoot{\@oddfoot}
\makeatother
\fi
\let\maketitlepage\makeagttitle
\let\makeshorttitle\maketitlepage
\let\maketitle\maketitlepage

\newwrite\gtoutfile
\long\gdef\makeheadfile{  
{\def\\{, }\def\s{ }
\immediate\openout\gtoutfile head.xxx
\immediate\write\gtoutfile{To: math@arxiv.org}
\immediate\write\gtoutfile{Subject: put OR rep NNNNN:ppppp}
\immediate\write\gtoutfile{--text follows this line--}
\immediate\write\gtoutfile{Proxy-for: \ifx\theasciiauthors\relax
\theauthors\else\theasciiauthors\fi\s<\ifx\theasciiemail\relax\theemail\else\theasciiemail\fi>}
\immediate\write\gtoutfile{\noexpand\\}
\immediate\write\gtoutfile{Authors: \ifx\theasciiauthors\relax
\theauthors\else\theasciiauthors\fi}
{\def\\{ }\immediate\write\gtoutfile{Title: \ifx\theasciititle\relax
\thetitle\else\theasciititle\fi}}
\immediate\write\gtoutfile{Subj-class: GT or SG, GR etc}
\immediate\write\gtoutfile{MSC-class: \theprimaryclass\ifx\thesecondaryclass\relax\else, \thesecondaryclass\fi}
\immediate\write\gtoutfile{Journal-ref: Algebr. Geom. Topol. \thevolumenumber\s
(\thevolumeyear) \startpage-\finishpage}
\immediate\write\gtoutfile{Comments: Published by Algebraic and
Geometric Topology at}
\immediate\write\gtoutfile{\s\s\s  http://www.maths.warwick.ac.uk/agt/AGTVol\thevolumenumber/agt-\thevolumenumber-\thepapernumber.abs.html}
\immediate\write\gtoutfile{\noexpand\\}
\immediate\write\gtoutfile{}
\ifx\theasciiabstract\relax
\immediate\write\gtoutfile{\theabstract}\else
\immediate\write\gtoutfile{\theasciiabstract}\fi
\immediate\write\gtoutfile{}
\immediate\write\gtoutfile{\noexpand\\}
\immediate\write\gtoutfile{}
\immediate\closeout\gtoutfile}}  

\def\maketitlepage{\makeagttitle\makeheadfile}
\let\makeshorttitle\maketitlepage
\let\maketitle\maketitlepage

\def\ifplaintex{\expandafter\ifx\csname documentclass\endcsname\relax}

\def\gtp{{\mathsurround=0pt\it $\cal G\mskip-2mu$eometry \&\ 
$\cal T\!\!$opology $\cal P\!$ublications}}  

\def\recd{{\small Received:\qua\receiveddate\ifx\reviseddate\relax
\else\qquad Revised:\qua\reviseddate\fi\par}} 

\def\lognumber#1{\def\thelognumber{#1}}
\def\volumenumber#1{\def\thevolumenumber{#1}}
\def\volumeyear#1{\def\thevolumeyear{#1}}
\def\papernumber#1{\def\thepapernumber{#1}}
\def\pagenumbers#1#2{\def\startpage{#1}\def\finishpage{#2}}
\def\published#1{\def\publishdate{#1}}

\def\received#1{\def\receiveddate{#1}}

\def\accepted#1{\def\accepteddate{#1}}

\def\asciiaddress#1{\def\theasciiaddress{#1}}

\long\def\asciiabstract#1{\long\def\theasciiabstract{#1}}

\let\\\par\let\thelognumber\relax\let\thevolumenumber\relax
\let\thepapernumber\relax\let\thevolumeyear\relax\let\startpage\relax
\let\finishpage\relax\let\publishdate\relax\let\receiveddate\relax
\let\reviseddate\relax\let\accepteddate\relax\let\theasciititle\relax
\let\theasciiauthors\relax\let\theasciiaddress\relax
\let\theasciiabstract\relax

\let\theasciiemail\relax

\ifplaintex
\font\logobig=cmssbx10 scaled 3836
\font\logomed=cmssbx10 scaled 2557
\else
\font\logobig=cmssbx10 scaled 4200
\font\logomed=cmssbx10 scaled 2800
\fi

\long\def\makeagttitle{   
\count0=\startpage
\agt\hfill      
\hbox to 45truept{\vbox to 0pt{\vglue -13truept{\logomed A\kern -.37em{\logobig 
T}\kern -.38em G}\vss}\hss}
\break
{\small Volume \thevolumenumber\ (\thevolumeyear)
\startpage--\finishpage\nl
Published: \publishdate}

\vglue .25truein

{\parskip=0pt\leftskip 0pt plus
1fil\def\\{\par\smallskip}{\Large\bf\thetitle}\par\medskip} \vglue
0.05truein

{\parskip=0pt\leftskip 0pt plus 1fil\def\\{\par}{\sc\theauthors}
\par\medskip}%
 
\vglue 0.03truein 

{\small\leftskip 25truept\rightskip 25truept{\bf Abstract}\stdspace\theabstract

{\bf AMS Classification}\stdspace\theprimaryclass
\ifx\thesecondaryclass\relax\else; \thesecondaryclass\fi\par
{\bf Keywords}\stdspace \thekeywords\par}\vglue 7truept

}   

\ifplaintex
\font\phead=cmsl9 scaled 950
\font\pnum=cmbx10 scaled 913
\font\pfoot=cmsl9 scaled 950
\headline{\vbox to 0pt{\vskip -4.5mm\line{\small\phead\ifnum
\count0=\startpage ISSN 1472-2739 (on-line) 1472-2747 (printed)
\hfill {\pnum\folio}\else\ifodd\count0\def\\{ }%
\ifx\theshorttitle\relax\thetitle\else\theshorttitle\fi\hfill{\pnum\folio}
\else\def\\{ and }{\pnum\folio}\hfill\ifx\theshortauthors\relax\theauthors
\else\theshortauthors\fi\fi\fi}\vss}}
\footline{\vbox to 0pt{\vglue 0mm\line{\small\pfoot\ifnum\count0=\startpage
\copyright\ \gtp\hfill\else
\agt, Volume \thevolumenumber\ (\thevolumeyear)\hfill\fi}\vss}}
\else
\font=cmsl9 scaled 1050
\font\lnum=cmbx10 
\font=cmsl9 scaled 1050
\makeatletter
\def\@oddhead{{\small\lhead\ifnum\count0=\startpage ISSN 1472-2739 
(on-line) 1472-2747 (printed)\hfill {\lnum\number\count0}\else\ifodd\count0
\def\\{ }\ifx\theshorttitle\relax \thetitle \else\theshorttitle\fi\hfill
{\lnum\number\count0}\else\def\\{ and }{\lnum\number\count0}
\hfill\ifx\theshortauthors\relax 
\theauthors\else\theshortauthors\fi\fi\fi}}\def\@evenhead{\@oddhead}
\def\@oddfoot{\small\lfoot\ifnum\count0=\startpage\copyright\ \gtp\hfill\else
\agt, Volume \thevolumenumber\ (\thevolumeyear)\hfill\fi}
\def\@evenfoot{\@oddfoot}
\makeatother
\fi
\let\maketitlepage\makeagttitle
\let\makeshorttitle\maketitlepage
\let\maketitle\maketitlepage

\newwrite\gtoutfile
\long\gdef\makeheadfile{  
{\def\\{, }\def\s{ }
\immediate\openout\gtoutfile head.xxx
\immediate\write\gtoutfile{To: math@arxiv.org}
\immediate\write\gtoutfile{Subject: put OR rep NNNNN:ppppp}
\immediate\write\gtoutfile{--text follows this line--}
\immediate\write\gtoutfile{Proxy-for: \ifx\theasciiauthors\relax
\theauthors\else\theasciiauthors\fi\s<\ifx\theasciiemail\relax\theemail\else\theasciiemail\fi>}
\immediate\write\gtoutfile{\noexpand\\}
\immediate\write\gtoutfile{Authors: \ifx\theasciiauthors\relax
\theauthors\else\theasciiauthors\fi}
{\def\\{ }\immediate\write\gtoutfile{Title: \ifx\theasciititle\relax
\thetitle\else\theasciititle\fi}}
\immediate\write\gtoutfile{Subj-class: GT or SG, GR etc}
\immediate\write\gtoutfile{MSC-class: \theprimaryclass\ifx\thesecondaryclass\relax\else, \thesecondaryclass\fi}
\immediate\write\gtoutfile{Journal-ref: Algebr. Geom. Topol. \thevolumenumber\s
(\thevolumeyear) \startpage-\finishpage}
\immediate\write\gtoutfile{Comments: Published by Algebraic and
Geometric Topology at}
\immediate\write\gtoutfile{\s\s\s  http://www.maths.warwick.ac.uk/agt/AGTVol\thevolumenumber/agt-\thevolumenumber-\thepapernumber.abs.html}
\immediate\write\gtoutfile{\noexpand\\}
\immediate\write\gtoutfile{}
\ifx\theasciiabstract\relax
\immediate\write\gtoutfile{\theabstract}\else
\immediate\write\gtoutfile{\theasciiabstract}\fi
\immediate\write\gtoutfile{}
\immediate\write\gtoutfile{\noexpand\\}
\immediate\write\gtoutfile{}
\immediate\closeout\gtoutfile}}  

\def\maketitlepage{\makeagttitle\makeheadfile}
\let\makeshorttitle\maketitlepage
\let\maketitle\maketitlepage

\def\ifplaintex{\expandafter\ifx\csname documentclass\endcsname\relax}

\def\gtp{{\mathsurround=0pt\it $\cal G\mskip-2mu$eometry \&\ 
$\cal T\!\!$opology $\cal P\!$ublications}}  

\def\recd{{\small Received:\qua\receiveddate\ifx\reviseddate\relax
\else\qquad Revised:\qua\reviseddate\fi\par}} 

\def\lognumber#1{\def\thelognumber{#1}}
\def\volumenumber#1{\def\thevolumenumber{#1}}
\def\volumeyear#1{\def\thevolumeyear{#1}}
\def\papernumber#1{\def\thepapernumber{#1}}
\def\pagenumbers#1#2{\def\startpage{#1}\def\finishpage{#2}}
\def\published#1{\def\publishdate{#1}}

\def\received#1{\def\receiveddate{#1}}

\def\accepted#1{\def\accepteddate{#1}}

\def\asciiaddress#1{\def\theasciiaddress{#1}}

\long\def\asciiabstract#1{\long\def\theasciiabstract{#1}}

\let\\\par\let\thelognumber\relax\let\thevolumenumber\relax
\let\thepapernumber\relax\let\thevolumeyear\relax\let\startpage\relax
\let\finishpage\relax\let\publishdate\relax\let\receiveddate\relax
\let\reviseddate\relax\let\accepteddate\relax\let\theasciititle\relax
\let\theasciiauthors\relax\let\theasciiaddress\relax
\let\theasciiabstract\relax

\let\theasciiemail\relax

\ifplaintex
\font\logobig=cmssbx10 scaled 3836
\font\logomed=cmssbx10 scaled 2557
\else
\font\logobig=cmssbx10 scaled 4200
\font\logomed=cmssbx10 scaled 2800
\fi

\long\def\makeagttitle{   
\count0=\startpage
\agt\hfill      
\hbox to 45truept{\vbox to 0pt{\vglue -13truept{\logomed A\kern -.37em{\logobig 
T}\kern -.38em G}\vss}\hss}
\break
{\small Volume \thevolumenumber\ (\thevolumeyear)
\startpage--\finishpage\nl
Published: \publishdate}

\vglue .25truein

{\parskip=0pt\leftskip 0pt plus
1fil\def\\{\par\smallskip}{\Large\bf\thetitle}\par\medskip} \vglue
0.05truein

{\parskip=0pt\leftskip 0pt plus 1fil\def\\{\par}{\sc\theauthors}
\par\medskip}%
 
\vglue 0.03truein 

{\small\leftskip 25truept\rightskip 25truept{\bf Abstract}\stdspace\theabstract

{\bf AMS Classification}\stdspace\theprimaryclass
\ifx\thesecondaryclass\relax\else; \thesecondaryclass\fi\par
{\bf Keywords}\stdspace \thekeywords\par}\vglue 7truept

}   

\ifplaintex
\font\phead=cmsl9 scaled 950
\font\pnum=cmbx10 scaled 913
\font\pfoot=cmsl9 scaled 950
\headline{\vbox to 0pt{\vskip -4.5mm\line{\small\phead\ifnum
\count0=\startpage ISSN 1472-2739 (on-line) 1472-2747 (printed)
\hfill {\pnum\folio}\else\ifodd\count0\def\\{ }%
\ifx\theshorttitle\relax\thetitle\else\theshorttitle\fi\hfill{\pnum\folio}
\else\def\\{ and }{\pnum\folio}\hfill\ifx\theshortauthors\relax\theauthors
\else\theshortauthors\fi\fi\fi}\vss}}
\footline{\vbox to 0pt{\vglue 0mm\line{\small\pfoot\ifnum\count0=\startpage
\copyright\ \gtp\hfill\else
\agt, Volume \thevolumenumber\ (\thevolumeyear)\hfill\fi}\vss}}
\else
\font=cmsl9 scaled 1050
\font\lnum=cmbx10 
\font=cmsl9 scaled 1050
\makeatletter
\def\@oddhead{{\small\lhead\ifnum\count0=\startpage ISSN 1472-2739 
(on-line) 1472-2747 (printed)\hfill {\lnum\number\count0}\else\ifodd\count0
\def\\{ }\ifx\theshorttitle\relax \thetitle \else\theshorttitle\fi\hfill
{\lnum\number\count0}\else\def\\{ and }{\lnum\number\count0}
\hfill\ifx\theshortauthors\relax 
\theauthors\else\theshortauthors\fi\fi\fi}}\def\@evenhead{\@oddhead}
\def\@oddfoot{\small\lfoot\ifnum\count0=\startpage\copyright\ \gtp\hfill\else
\agt, Volume \thevolumenumber\ (\thevolumeyear)\hfill\fi}
\def\@evenfoot{\@oddfoot}
\makeatother
\fi
\let\maketitlepage\makeagttitle
\let\makeshorttitle\maketitlepage
\let\maketitle\maketitlepage

\newwrite\gtoutfile
\long\gdef\makeheadfile{  
{\def\\{, }\def\s{ }
\immediate\openout\gtoutfile head.xxx
\immediate\write\gtoutfile{To: math@arxiv.org}
\immediate\write\gtoutfile{Subject: put OR rep NNNNN:ppppp}
\immediate\write\gtoutfile{--text follows this line--}
\immediate\write\gtoutfile{Proxy-for: \ifx\theasciiauthors\relax
\theauthors\else\theasciiauthors\fi\s<\ifx\theasciiemail\relax\theemail\else\theasciiemail\fi>}
\immediate\write\gtoutfile{\noexpand\\}
\immediate\write\gtoutfile{Authors: \ifx\theasciiauthors\relax
\theauthors\else\theasciiauthors\fi}
{\def\\{ }\immediate\write\gtoutfile{Title: \ifx\theasciititle\relax
\thetitle\else\theasciititle\fi}}
\immediate\write\gtoutfile{Subj-class: GT or SG, GR etc}
\immediate\write\gtoutfile{MSC-class: \theprimaryclass\ifx\thesecondaryclass\relax\else, \thesecondaryclass\fi}
\immediate\write\gtoutfile{Journal-ref: Algebr. Geom. Topol. \thevolumenumber\s
(\thevolumeyear) \startpage-\finishpage}
\immediate\write\gtoutfile{Comments: Published by Algebraic and
Geometric Topology at}
\immediate\write\gtoutfile{\s\s\s  http://www.maths.warwick.ac.uk/agt/AGTVol\thevolumenumber/agt-\thevolumenumber-\thepapernumber.abs.html}
\immediate\write\gtoutfile{\noexpand\\}
\immediate\write\gtoutfile{}
\ifx\theasciiabstract\relax
\immediate\write\gtoutfile{\theabstract}\else
\immediate\write\gtoutfile{\theasciiabstract}\fi
\immediate\write\gtoutfile{}
\immediate\write\gtoutfile{\noexpand\\}
\immediate\write\gtoutfile{}
\immediate\closeout\gtoutfile}}  

\def\maketitlepage{\makeagttitle\makeheadfile}
\let\makeshorttitle\maketitlepage
\let\maketitle\maketitlepage

\lognumber{39}
\volumenumber{2}
\volumeyear{2002}
\papernumber{39}
\pagenumbers{949}{1000}
\received{2 August 2002}
\accepted{12 October 2002}
\published{25 October 2002}

\usepackage{amsmath}  
\usepackage{amssymb}
\usepackage{graphics}
\usepackage[all]{xy}

\numberwithin{equation}{section}

\def\bea{\begin{align}}
\def\eea{\end{align}}

\newtheorem{Thm}{Theorem}[section]
\newtheorem*{Thm*}{Theorem}
\newtheorem{Prop}[Thm]{Proposition}
\newtheorem{Lem}[Thm]{Lemma}
\newtheorem{Cor}[Thm]{Corollary}

\theoremstyle{remark}
\newtheorem{Rem}[Thm]{Remark}
\newtheorem*{Rem*}{Remark}
\newtheorem*{Ack}{Acknowledgments}

\theoremstyle{definition}
\newtheorem{Def}[Thm]{Definition}

\newtheorem{Conj}[Thm]{Conjecture}
\newtheorem*{Convs}{Conventions}
\newtheorem{Conv}[Thm]{Convention}

\def\bconv{\begin{Conv}}
\def\econv{\end{Conv}}
\def\bconvs{\begin{Convs}}
\def\econvs{\end{Convs}}
\def\bprop{\begin{Prop}}
\def\eprop{\end{Prop}}
\def\bconj{\begin{Conj}}
\def\econj{\end{Conj}}
\def\bth{\begin{Thm}}
\def\eth{\end{Thm}}
\def\bcor{\begin{Cor}}
\def\ecor{\end{Cor}}
\def\blem{\begin{Lem}}
\def\elem{\end{Lem}}
\def\bdefiniz{\begin{Def}}
\def\edefiniz{\end{Def}}
\def\incul{\hookrightarrow}

\DeclareMathOperator{\diff}{Diff}

\DeclareMathOperator{\ord}{ord}

\newcommand{\lora}{{\longrightarrow}}

\newcommand{\rsa}{\mapsto}

\newcommand\de{\partial}
\newcommand{\R}{\mathbb{R}}

\newcommand{\N}{\mathbb{N}}
\newcommand{\reali}{\mathbb{R}}
\newcommand{\bbR}{\mathbb{R}}

\newcommand{\naturali}{\mathbb{N}}

\newcommand{\imb}[1]{{\,\hbox{\rm Imb}\,(S^1,{#1}) }}
\newcommand{\imm}[1]{{\,\hbox{\rm Imm}\,(S^1,{#1}) }}

\newcommand{\imbr}[1]{{\,\hbox{\rm Imb}\,(S^1,\reali^{#1}) }}
\newcommand{\imbrf}[1]{{\,\hbox{\rm Imb}_\mathrm{f}\,(S^1,\reali^{#1}) }}
\newcommand{\immr}[1]{{\,\hbox{\rm Imm}\,(S^1,\reali^{#1}) }}
\newcommand{\immrf}[1]{{\,\hbox{\rm Imm}_\mathrm{f}\,(S^1,\reali^{#1}) }}

\newcommand{\immsr}[2]{{\,\hbox{${\rm Imm}_{#1}$}\,(S^1,\reali^{#2}) }}
\newcommand{\immpsr}[2]{{\,\hbox{${\rm Imm}_{#1}^{\prime}$}\,(S^1,\reali^{#2}) }}

\def\cmp#1,{{ Commun.\ Math.\ Phys.\ \bf #1},}
\def\jmp#1,{{ J.\ Math.\ Phys.\ \bf #1},}
\def\pl#1,{{ Phys.\ Lett.\ \bf #1},}
\def\npb#1,{{ Nucl.\ Phys.\ {\bf B #1}},}
\def\mpl#1,{{ Mod.\ Phys.\ Lett.\ \bf #1},}
\def\pr#1,{{ Phys.\ Rev.\ \bf #1},}
\def\prl#1,{{ Phys.\ Rev.\ Lett.\ \bf #1},}
\def\lmp#1,{{ Lett.\ Math.\ Phys.\ \bf #1},}
\def\jktr#1,{{ J.\ of Knot Theory and its Ramifications\ \bf #1},}
\def\bams#1,{{ Bull.\ Amer.\ Math.\ Soc.\ \bf #1},}
\def\jdg#1,{{ J. Diff.\ Geom.\ \bf #1},}
\def\asm#1,{{Adv.\ Sov.\ Math.\ \bf #1},}
\def\gt#1,{{Geometry\ \&\ Topology \bf #1},}

\begin{document}

\title{Configuration spaces and Vassiliev classes\\in any dimension}

\author{Alberto S. Cattaneo\\Paolo Cotta-Ramusino\\Riccardo Longoni}
\shortauthors{Cattaneo, Cotta-Ramusino and Longoni}

\addresses{Mathematisches Institut, Universit\"at Z\"urich--Irchel,
Winterthurerstrasse 190\\CH-8057 Z\"urich, Switzerland\\\smallskip
\\Dipartimento di 
Fisica, Universit\`a degli Studi di
Milano \& INFN Sezione di Milano\\Via Celoria, 16, I-20133 
Milano, Italy\\\smallskip
\\Dipartimento di Matematica ``G. Castelnuovo'',
Universit\`a di Roma ``La Sapienza''\\Piazzale Aldo Moro, 5,
I-00185 Roma, Italy} 

\asciiaddress{Mathematisches Institut, Universitat Zurich-Irchel,
Winterthurerstrasse 190\\CH-8057 Zurich, Switzerland\\Dipartimento di 
Fisica, Universita degli Studi di
Milano and INFN Sezione di Milano\\Via Celoria, 16, I-20133 
Milano, Italy\\Dipartimento di Matematica `G. Castelnuovo',
Universita di Roma `La Sapienza'\\Piazzale Aldo Moro, 5,
I-00185 Roma, Italy}

\email{asc@math.unizh.ch, cotta@mi.infn.it, longoni@mat.uniroma1.it}

\begin{abstract}
The real cohomology of the space of imbeddings of $S^1$ into
$\reali^n$, $n>3$, is studied by using configuration space integrals.
Nontrivial classes are explicitly constructed. As a by-product, we
prove the nontriviality of certain cycles of imbeddings obtained by
blowing up transversal double points in immersions. These cohomology
classes generalize in a nontrivial way the Vassiliev knot invariants.
Other nontrivial classes are constructed by considering the
restriction of classes defined on the corresponding spaces of
immersions.
\end{abstract}

\asciiabstract{
The real cohomology of the space of imbeddings of S^1 into
R^n, n>3, is studied by using configuration space integrals.
Nontrivial classes are explicitly constructed. As a by-product, we
prove the nontriviality of certain cycles of imbeddings obtained by
blowing up transversal double points in immersions. These cohomology
classes generalize in a nontrivial way the Vassiliev knot invariants.
Other nontrivial classes are constructed by considering the
restriction of classes defined on the corresponding spaces of
immersions.}



\maketitle

\section{Introduction}

In this paper we study de~Rham cohomology classes of the space $\imbr n$
of smooth imbeddings of $S^1$ into $\bbR^n$, $n>3$, using as main tools
configuration-space integrals and graph cohomology.

Before describing the setting of this paper we give a brief description
of the main results obtained.

\subsection{Main results}
We consider two complexes  $(\mathcal D_o^{k,m}, \delta_o)$ and
$(\mathcal D_e^{k,m},\delta_e)$ generated by some decorated
graphs. These {\em graph complexes}\/ are bigraded by two integers
$m$, $k$ called respectively the degree and the order. The order
is minus the Euler characteristic of the graph, while the degree
measures the deviation of the graph from being trivalent. The
coboundary operators increase the degree by one and do not change
the order.

We prove in subsection~\ref{finale}
the following:
\begin{Thm}\label{uno}
For every $k\in \N$, there exist chain maps from graph complexes to
de~Rham complexes
\begin{equation}
\label{ch-map-intro-o}
\mathcal D^{k,m}_o \to \Omega^{(n-3)k+m}(\imbr n)\qquad\mbox{for $n$ odd}
\end{equation}
\begin{equation}
\label{ch-map-intro-e}
\mathcal D^{k,m}_e \to \Omega^{(n-3)k+m}(\imbr n)\qquad\mbox{for $n$ even}
\end{equation}
that induce injective maps in cohomology when $m=0$.
\end{Thm}
From the combinatorial structure of the graph complexes, one immediately
deduces (see again subsection~\ref{finale}) the following:
\begin{Cor}
\label{cor-uno} For any $n> 3$ and for any positive integer $k_0$,
there are nontrivial cohomology classes on $\imbr{n}$
of degree greater than $k_0$.
\end{Cor}
All the differential forms appearing in the above
Theorem turn out to be equivariant w.r.t.\ the action of the group
$\mathit{Diff}^+(S^1)$ of orientation preserving diffeomorphisms
of the circle.

The classes of Theorem~\ref{uno} can be seen as an extension to
higher dimensions of Vassiliev knot invariants. One of the main
ingredients of Vassiliev's approach is to consider immersions
which are imbeddings but for a finite number (say $k$) of
transversal double points (let us call them ``special
immersions''). The ways of pushing off a double point form, up to
homotopy, an $(n-3)$-dimensional sphere (viz., one must choose a
normalized vector transverse to the plane spanned by the tangent
vectors of the intersecting strands at the double point). So every
special immersion with $k$ double points gives rise to a
$k(n-3)$-cycle in $\imbr n$. Our construction allows us to prove
that infinitely many of these cycles are nontrivial.

When we extend the above construction to imbeddings  with framing,
namely when we consider pairs $(K,w)$ consisting of an imbedding
$K\colon S^1\to \bbR^n$ and a section $w$ of the pulled-back bundle
$K^*SO(\bbR^n)$, then the situation becomes simpler and we prove in
subsection~\ref{prooftre} the following:
\begin{Thm}\label{tre}
All cycles of framed imbeddings of $S^1$ into $\bbR^{2s+1}$
determined by framed special immersions are nontrivial.
\end{Thm}

The restriction of cohomology classes on the space of immersions
$\immr n$ to the space of imbeddings $\imbr n$ is also discussed.
Combining Theorems~\ref{nodd}, Proposition~\ref{p1veven} and
Corollary~\ref{p1-vn-even} we obtain:
\begin{Thm}\label{due}
When $n$ is odd, the inclusion $\imbr n \incul \immr n$
induces the zero map in cohomology. When $n$ is even, the inclusion
map $\imbr {n}\incul\immr {n}$ is nontrivial in cohomology.
\end{Thm}
Contrary to what happens in Thm.~\ref{uno}, not all the
differential forms of Thm.~\ref{due} are
$\mathit{Diff}^+(S^1)$-equivariant. However, the equivariant ones
turn out to be in the image of certain graphs (of degree different
from zero) through the map \eqref{ch-map-intro-e}. Thus,
Thm.~\ref{due} provides an extension of the last statement of
Thm.~\ref{uno} in the case $m\neq 0$.

\subsection{The setting}
The configuration spaces $C_q^0(M)$ of a manifold $M$ are simply
the Cartesian products $M^q$ minus all diagonals.
Configuration spaces are naturally associated to imbeddings. Indeed, if
$f\colon N\to M$ is an imbedding, it defines maps
$C_q^0(f)\colon C_q^0(N)\to C_q^0(M)$ for every $q\ge0$.

This simple, natural relation has an application to knot
invariants, i.e., to the study of the zeroth cohomology of the
space of imbeddings of $S^1$ into $\bbR^3$. We refer to Bott and
Taubes's \cite{BT} construction inspired from the perturbative
expansion of Chern--Simons theory \cite{Wit}. The ``physical''
origin of this construction should not be too surprising; for, as
a matter of fact, a ``correlation function'' in physics (i.e., the
inverse of some differential operator) is usually well-defined
only at non-coincident points---and this also leads naturally to
configuration spaces.

The gist of their construction is as follows: One considers differential
forms on $C_q^0(\bbR^3)$ given by products of the
rotational-invariant representatives
of the two-dimensional
generators in cohomology (``tautological forms'');
next one integrates these forms on cycles defined by constraining
some of the $q$ points to lie in the image of the given imbedding $K$;
finally, one wants to
prove that certain linear combinations of these integrals
are actually invariant under isotopies of $K$.

The big technical problem in Bott and Taubes's construction (as
well as in other related constructions \cite{AS,K,BC1,BC2} for
$3$-manifold invariants) is the convergence of the above
integrals, a nontrivial fact since the tautological forms are not
compactly supported. An elegant solution---which also lies at the
core of the subsequent construction of invariants by determining
the suitable linear combinations of integrals---relies on a
compactification of configuration spaces on which the tautological
forms extend as smooth forms. This is the Fulton--MacPherson
\cite{FMP} compactification, but in the differential-geometric
version later given by Axelrod and Singer \cite{AS}. What is
actually needed is still a further improvement; viz., a
compactification of the configuration spaces of $\bbR^3$ with some
points lying on the knot. This is done in \cite{BT}, and indeed in
the more general case of the configuration space of a manifold $M$
with some points lying on the image of a given imbedding of
another manifold $N$ (assuming $N$ and $M$ to be compact). The
last result allows one to approach more general imbedding
problems, as is done in the present work. We finally observe that,
in the $3$-dimensional case, only knot invariants have been
constructed this way (and no higher-degree cohomology classes on
the space of imbeddings) and, moreover, that these invariants are
proved to exist and to be nontrivial, but are obtained modulo
corrections with unknown coefficients.

Let us now turn to the case of imbeddings of $S^1$ into $\bbR^n$
with $n>3$. Here ``tautological forms'' are representatives of the
$(n-1)$-dimensional cohomology generators of configuration spaces
of $\bbR^n$. We integrate products of tautological forms on cycles
in configuration spaces constraining some of the points to lie on
the imbedding, thus getting differential forms on $\imbr n$. We
construct closed linear combinations of these forms by using graph
cohomology as explained in the following.

Graphs, with edges corresponding to tautological forms, are a
simple way of keeping track of all configuration-space integrals
one may consider. In order to get a complete description, i.e.\
without sign ambiguities, one actually has to decorate the graphs
in a certain way (in fact, two different ways corresponding to $n$
even or odd). At this point, one can define a grading and a
coboundary operator on the vector space generated by all decorated
graphs in such a way that the assignment to each graph of the
corresponding configuration-space integral defines a (degree
shifting) chain map from the graph complex to de~Rham complex of
$\imbr n$, see Theorem~\ref{thm-Ivn}. We give an explicit
definition of graph cohomologies, along the lines originally
proposed by Kontsevich \cite{K}, and describe in details the chain
maps and so the relation to the problem of imbeddings. It is at
this point that it is crucial to have $n>3$, as for $n=3$ the map
we construct might fail to be a chain map. As was
observed by the referee, this chain map is also a morphism of
differential graded commutative algebras, with the multiplication
on the graphs defined as (a graded version of) the shuffle
product. We plan to return on this and to discuss other algebraic
structures on the graph complex in \cite{CCRL2}. Finally, we
prove that the induced map in cohomology is injective in degree
zero. This is done by pairing the corresponding cohomology classes
on $\imbr n$ to the cycles arising from special immersions
described before. Observe that, for $n>3$, the connected
components of the space of special immersions are in one-to-one
correspondence with chord diagrams (i.e., graphs consisting of a
distinguished circle with chords), each chord representing a
transversal double point. We prove that the pairing of a cycle of
imbeddings determined by a special immersion with a differential
form coming from a graph cocycle containing the corresponding
chord diagram is non zero, see Theorem~\ref{HD-HImm}. Moreover, we
prove that every graph cocycle of degree zero contains a chord
diagram, see Proposition~\ref{thm-trichord}.

Consider now the space $\imbrf n$ of framed imbeddings. Since we have a
projection to $\imbr n$ (forgetting the framing), we can pull back
all cohomology classes constructed before. Given a framed special
immersion, we can then generate a cycle of framed imbeddings and exactly
as above prove that cohomology classes corresponding to graph cocycles
in degree zero are nontrivial.
In odd dimensions,
we can extend our construction and produce new classes (corresponding
to new classes in a modified graph cohomology), and again prove nontriviality
in degree zero. But at this point, we can also use the same technique
to prove that all cycles determined by framed special immersions are
nontrivial, as stated in the Theorem at the beginning.

A problem which is strictly related to the subject of this paper
is the study of the cohomology of the spaces of {\em long knots},
i.e., the spaces of imbeddings of the real line into Euclidean
spaces with fixed behavior at infinity. This problem has been
already addressed by various authors \cite{GW, Si, T, V3}. We
comment here that the methods developed in this paper are easily
generalized in that direction. In particular one can prove
\cite{Tpriv} that our graph complexes are quasi-isomorphic to the
first term of the spectral sequences defined by Vassiliev \cite{T,
V3} and Sinha \cite{Si}, and that Theorem~\ref{uno} implies the
convergence of these spectral sequences along the diagonal.

We conclude with some remarks on the relation between the
configuration\hspace{0pt}-\hspace{0pt}space techniques described
above and physics. We recall that in the $3$-dimensional case,
Vassiliev knot invariants appear in the perturbative expansion of
expectation values of traces of holonomies in Chern--Simons theory
\cite{Wit}. Bott and Taubes's construction is based on the
expansion in the ``covariant gauge'' \cite{GMM,BN90}, whereas the
Kontsevich integral \cite{KVass,BN95} is based on the expansion in
the ``holomorphic gauge'' \cite{FK}. As described in
\cite{CCRFM95}, the same Vassiliev invariants may also be obtained
in the perturbative expansion of $BF$ theory in $3$ dimensions.
This theory can actually be defined in any dimension, and the
perturbative expansion of expectation values of traces of the
generalized holonomies defined in \cite{CCRRo,CRo} (see also
\cite{CCRR}) are actually related to the cohomology classes of
imbeddings described in the present paper. Moreover, the analysis
of the $BF$ theories made in \cite{CFP} suggests the possibility
of connecting our results with the string topology of Chas and
Sullivan \cite{CS}.

\subsection{Plan of the paper}

In Section~\ref{vorder} we define a map that assigns to each element
of $H^p(\imbr n, \reali)$ a cohomology class in
$H^{p-k(n-3)}(\immsr k n,\reali)$ where $\immsr k n$ denotes the
space of immersion with exactly $k$ transversal double points.
This map is the generalization of the map that extends knot
invariants to invariants of $\immsr k 3$ \cite{V1,BN95}. We say
then that a real cohomology class of $\imbr n$ has  {\em
Vassiliev-order} $s$ if the corresponding cohomology class of
$\immsr k  n$ is zero for $k>s$ and non-zero for $k=s$.

After recalling the Bott--Taubes
construction for tautological forms and configuration spaces in
Section~\ref{bott-taubes},
we define in Section~\ref{ddco} the two complexes
$(\mathcal D_o,\delta_o)$ and $(\mathcal D_e,\delta_e)$
mentioned above.
We show that the configuration space integral is
a chain map from the above complexes to the
de~Rham complex of $\imbr n$ (where the two cases $n$ even and $n$ odd
are kept separately).

In Section~\ref{sect-triv} we focus on  trivalent graphs and
construct explicitly some nontrivial cocycles  that are given by
linear combinations of them.

In Section~\ref{nontriv} we show that for $n>3$ the morphisms between
the complexes $(\mathcal D_o,\delta_o)$, $(\mathcal D_e,\delta_e)$
and $(\Omega^*(\imbr n), d)$ ($n$ odd and, respectively, even) are
monomorphisms in cohomology when they are restricted to the subspaces of
trivalent graphs. At the end, we prove Thm.~\ref{uno} and Cor.~\ref{cor-uno}.

In Section~\ref{sec-framed} we consider the space of framed imbeddings
$\imbrf n$. We show how to define new classes in case $n$ is odd.
Here we define a modified graph cohomology and a new chain map, which
we prove to be injective in cohomology in degree zero.
We also prove Thm.~\ref{tre}.

In Section~\ref{chen} we recall the construction of the
generators of $H^*(\immr n, \reali)$ via Chen integrals
\cite{Chen} and compute their restrictions to $\imbr n$. We show
that this restriction is trivial if $n$ is odd but yields
nontrivial classes of $\imbr n$ if $n$ even, thus proving Thm.~\ref{due}.

Finally, in the Appendix we discuss some Vanishing Theorems that are
needed in order to define  the morphisms between the complexes considered
before. The main result is that, in computing the differential of an
integral of tautological forms, contributions from the so-called
hidden faces are always zero.\\

\bconvs
Throughout this paper we assume $n> 3$, unless otherwise stated.

We also assume that all the spaces under consideration (namely,
$S^1$ and $\reali^n$) are oriented. In particular two imbeddings
(or immersions) that are obtained from each other by reversing the
orientation of $S^1$ will be considered as different
elements of $\imbr n$ (or $\immr n$).

We are concerned only with {\em real}\/ cohomology groups that we will
denote by
$H^*(\imbr n)$ or $H^*(\immr n)$.

Finally, in the course of the paper we need to choose a unit
generator of the top cohomology of $S^n$. The main results of
Section~\ref{chen} are independent of such choice. In the rest of
this paper, however, we need to restrict ourselves to {\em
symmetric forms}\/ (see Definition~\ref{def-gf}). \econvs

\begin{Ack}
We thank especially Raoul Bott, Jim Stasheff, Victor Vassiliev and
the referee for pointing out parts of the previous version of this
paper that needed clarification or corrections. We thank Giovanni
Felder, Nathan Habegger, Maurizio Rinaldi, Dev Sinha, Dennis
Sullivan and Victor Tourtchine for useful comments and interesting
discussions. We thank Carlo Rossi and Simone Mosconi for carefully
reading the manuscript.

A.~S.~C. thanks the I.N.F.N., Sezione di Milano,
P.~C.-R. thanks the Universit\"at Z\"urich--Irchel and the ETH Z\"urich,
and R.~L. thanks the Universit\"at Z\"urich--Irchel
for their hospitality.

A.~S.~C. thanks partial support of SNF Grant No.$\backslash
2100-055536.98/1$. P.~C.-R. and R.~L. thank partial support of
MURST.
\end{Ack}

\section{Vassiliev classes in $H^*(\imbr n)$}
\label{vorder}

In this Section we propose a classification scheme for the
cohomology classes in $H^*(\imbr n)$, including those that not are
necessarily obtained by pullback of classes in $H^*(\immr n)$ via
the inclusion map
\[
\imbr n\incul \immr n.
\]
This scheme is a direct
generalization of the scheme proposed by Vassiliev \cite{V1}
for knot invariants in $\reali^3$.

We consider the space
$\immsr k n$ which is defined as  the submanifold of $\immr n$
whose elements  have exactly $k$ transversal
double points. Moreover we set $\immpsr k n$ to be the submanifold
of $\immsr k n$
given by those immersions whose
initial point $\gamma(0)$ does not coincide with
any double point.

We enumerate all the double points of
any $\gamma\in \immpsr k n$ starting  from the initial point $\gamma(0)$.
Then we blow up, in order, all the
double points in the way  described below.

Let $\mathbf{x}^j=\gamma(t_1^j)=\gamma(t_2^j)$ be the $j$th double point,
with $t_1^j<t_2^j$. We denote by $l_1^j=D\gamma(t_1^j)$
and $l_2^j=D\gamma(t_2^j)$ the normalized tangent vectors
at $\mathbf{x}^j$ and by $T^j$ the plane in $T_{\mathbf{x}^j}\reali^n$
spanned by $l_1^j$ and $l_2^j$ with the orientation determined by
$l_1^j\wedge l_2^j$.

Here we assume to have chosen once and for all an orientation
in $\reali^n$. Moreover, for the rest of this section it is useful
to pick up a metric on $\reali^n$ as well.

Then we consider the $(n-2)$-plane $N^j\subset T_{\mathbf{x}^j}\reali^n$
that is normal to $T^j$ with the induced orientation and the space $Q^j$
of normalized vectors in $N^j$.

The space  ${\mathfrak Q}^j_k$ of pairs $(\gamma, \mathbf {z^j})$
of immersions with $k$ transversal double points
and normalized vectors in $Q^j$  can be formally
described as follows.
If we consider the Grassmann manifold $SG_{2,n}$
of {\em oriented}\/ 2-planes in $\reali^n$, then we have smooth maps
\[
r_{k,j}:\immpsr k n\to SG_{2,n}\equiv SO(n)/\{SO(2)\times SO(n-2)\}
\]
that associate to the $j$-th double point the oriented plane $T^j$.

We have an associated bundle
\[
Q=SO(n)\times_{\{SO(2)\times  SO(n-2)\}}
S^{n-3}\to SG_{2,n}
\]
whose fiber is the homogeneous space
\[
S^{n-3}=[SO(2)\times SO(n-2)]/[SO(2)\times SO(n-3)]\equiv SO(n-2)/SO(n-3).
\]
The  space
$Q$ can equivalently be
obtained by dividing $SO(n)$
by $SO(2)\times  SO(n-3)$.

The pull-back bundle ${\mathfrak Q}^j_k\equiv r_{k,j}^* Q$ is a sphere bundle
with fiber $S^{n-3}$
so that  the following diagram:
\[
\begin{array}{ccc}
{\mathfrak Q}^j_k&\lora & Q \\
\downarrow & & \downarrow\\
\immpsr k n & \stackrel{r_{k,j}}{\lora} & SG_{2,n} \\
\end{array}
\]
is commutative.

By considering in their order
all the double points, we can define the map
\[
r_k:\immpsr k n \to  \left(SG_{2,n}\right)^{\times k}
\]
and  the bundle
\begin{equation}
{\mathfrak Q}_k\equiv r_k^*Q^{\times k}
\label{frakq}
\end{equation}
with fiber $(S^{n-3})^{\times k}$.

Next we will define, see \eqref{zymap}, a map $s_k: {\mathfrak Q}_k\to
\imbr n$ that corresponds to the blow-up  of
all the double points of an immersion.

Given $\gamma\in \immpsr k n $ and  an element of the fiber
of  ${\mathfrak Q}^i_k$ over $\gamma$, which we represent as
$\mathbf{z}^j\in S^{n-3}$,
we  choose $a_j$
to be either $1$ or $2$ and define the following loops
in  $T_{\mathbf{x}^j}\reali^n$:
\begin{equation}
\alpha^{j}_{a_j}(\mathbf{z}^j)(t)=
\begin{cases}
0 & \text{if $t\notin [t^j_{a_j}-\epsilon,t^j_{a_j}+\epsilon]$},\\
(-1)^{a_j+1}\, \mathbf{z}^j\,\delta\exp
\left(1/[(t-t^j_{a_j})^2-\epsilon^2]\right)
&\text{if $t\in [t^j_{a_j}-\epsilon,t^j_{a_j}+\epsilon]$},
\end{cases}
\label{alfai}
\end{equation}
with $\epsilon,\delta >0$.

If we add to the immersion $\gamma$ the loop
$\alpha^{j}_{a_j}(\mathbf{z}^j)$,
using the natural identification $\reali^n\approx T_{\mathbf{x}^j}\reali^n$,
we remove the $j$th double point (see figure~\ref{figura1}).

We assume, from now on, that the parameters $\epsilon$ and $\delta$ are chosen
so small that no new double point is created by this operation.

In this construction one of the two strands that meet in the $j$th
double point is ``lifted'' in a way parameterized by $\mathbf{z}^j$
that belongs to the fiber over $\gamma$ of
the sphere bundle ${\mathfrak Q}^j_k$.
The union of all the possible lifts (for a given immersion $\gamma$
and a given
double point) describes the suspension
of the fiber $S^{n-3}$, namely, an $(n-2)$-sphere $\mathcal{S}^j_{a_j}$.
Denoting by $\ell^j_{a_j}$ the straight line passing through $\mathbf{x}^j$
with tangent $l^j_{a_j}$, we have the following
\bprop
The linking number between $\ell^j_{a_j}$ and\/ $\mathcal{S}^j_{b_j}$,
$b_j\equiv a_j+1\mod2$, is one.
\label{convicted}
\eprop
The proof is just a consequence of the orientation choices.
Observe that, being $\ell^j_{a_j}$ a 1-manifold, the above linking number
does not depend on the order.

For any given  choice of $a$ and of the ``small''
parameters $\epsilon$ and $\delta$ at each double point, we have thus
defined a map
\begin{equation}
s_k: {\mathfrak Q}_k\to \imbr n
\label{zymap}
\end{equation}
which is described, in any coordinate  neighborhood  of
$\gamma\in\immpsr k n$, by cycles:
\begin{equation}
(S^{n-3})^k\ni (\mathbf{z}^1,\cdots,\mathbf{z}^k)\rsa \gamma +
\sum_{j=1}^{k}\alpha^{j}_{a_j}(\mathbf{z}^j).
\label{zyklon}
\end{equation}
Due to  the arbitrariness in the choice of the index
$a_j\in\{1,2\}$ attached to each double point of $\gamma\in
\immpsr k n$, we have constructed $2^k$ cycles, for which we have
the following \bprop If we choose different values of
$a_j\in\{1,2\}$ for the double point labelled by $j$ in
(\ref{alfai}), then the resulting cycles (\ref{zyklon}) are
homologous. \label{homseg} \eprop
\begin{proof}
It is enough to consider two segments $[0,1]\ni t\rsa l_1^j(t)$
and $[0,1]\ni s\rsa l_2^j(s)$ that intersect
transversally at the middle point.  We choose $\mathbf{z}^j\in S^{n-3}$
and remove the crossing point as follows:
\begin{equation}
\left\{
\begin{array}{ccc}
l_1^j(t) & \rsa & l_1^j(t)+\mathbf{z}^j\delta\exp
\left(1/[(t-1/2)^2-\epsilon^2]\right)\\
l_2^j(s) & \rsa & l_2^j(s)- \mathbf{z}^j\delta\exp
\left(1/[(s-1/2)^2-\epsilon^2]\right)
\end{array}
\right.
\label{biciclo}
\end{equation}
where $\delta$ and $\epsilon$ are small positive numbers.
\begin{figure}[ht!]
\cl{\resizebox{8cm}{!}{\includegraphics{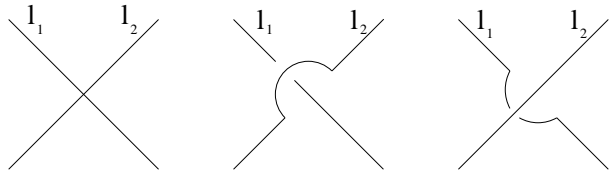}}}
\caption{The resolution of a transversal double intersection}
\label{figura1}
\end{figure}

We have then an $(n-3)$-cycle of imbedded pairs of segments with
fixed end-points. Let us now take $h\in [0,1]$. If we replace
$\delta$ with $h\delta$ in the first line of \eqref{biciclo}, then
we have a homotopy between the cycle (\ref{biciclo}) and the cycle
obtained by modifying only $l_2^j$. Analogously if we replace
$\delta$ with $h\delta$ in the second line of \eqref{biciclo}, we
have an homotopy between \eqref{biciclo} and the cycle obtained by
modifying only $l_1^j$.
\end{proof}

By pulling  back cohomology classes via  \eqref{zymap}
and integrating them along the fibers in ${\mathfrak Q}_k$
we  obtain the following morphisms in cohomology:
\begin{equation}
i_k^\prime:H^{p}(\imbr n)\to H^{p-k(n-3)}(\immpsr k n).
\label{zymor}
\end{equation}

Notice that the maps \eqref{zymor} are independent
of the choices of the $a$'s at each double point.

For future purposes, we extend the map
(\ref{zymor}) by setting it equal to zero if $p-k(n-3)<0$.
Hence $i_k^\prime$ is defined for every $k\in\naturali$.

\bdefiniz
We say that a cohomology class
$\omega\in H^p(\imbr n)$ is of  finite {\em Vassiliev-order}\/ (or
V-order) $k$ if  $i_s^\prime(\omega)= 0$ for every $s>k$ and
$i_k^\prime(\omega)\neq 0$. If $i_k^\prime(\omega)$ is non zero for
any $k$, then we say that the V-order is infinite.
\edefiniz

\begin{Rem}
If $n>3$, then the V-order is always finite. If $n=3$
then the V-order may be infinite.
The case $n=3$ and $p=0$ is the case of  knot invariants,
as originally studied by Vassiliev \cite{V1}.
\end{Rem}

If $n>3$ and $p=k(n-3)$, then
from (\ref{zymor})
we conclude that there is a morphism:
\begin{equation}
i_k^\prime :H^{k(n-3)}(\imbr n)\to H^0(\immpsr k n).
\label{hzero1}
\end{equation}
This case will be particularly important in the rest of the paper,
basically because of the following
\bprop
If $n>3$, then:
\begin{itemize}
\item[\rm(i)] the connected components of $\immsr k n$ are in one-to-one
correspondence with the set of chord diagrams with $k$ chords;
\item[\rm(ii)] the connected components of $\immpsr k n$ are in one-to-one
correspondence with chord diagrams with $k$ chords and a marked
point distinct from the end-points of the chords.
\end{itemize}
\label{hhzero} \eprop Here and in the following, by {\em chord
diagram}\/ we mean a circle with  chords that have no end-points
in common.

\begin{proof}
If $n>3$, then any finite
collection of (piecewise) imbedded  loops
can be isotopically deformed to a trivial link.
Hence the connected components of $\immsr k n$ are
determined uniquely by the position of the double points.
Their pre-images are points on a circle that are identified in pairs, i.e.,
chords.

In the case of $\immpsr k n$, one has just to take care of the additional
information given by the initial point.
\end{proof}

In general we want to determine whether a given class in $H^p(\imbr n)$
is trivial or not. The relevance of the order $p=k(n-3)$
is highlighted  by the following criterion:
\bcor
A sufficient condition for a class $\omega\in H^{k(n-3)}(\imbr n)$
to be nontrivial is that its image under (\ref{hzero1}) is nontrivial.
\label{crithzero}
\ecor

\begin{Rem}
\label{rem-eq}
We have a map
\begin{equation}
\label{def-varphi}
\varphi:H_0(\immsr k n)\to H_0(\immpsr k n)
\end{equation}
which associates to any chord diagram $D$ 
the average  of all inequivalent
chord diagrams with a marked point  that have the same chords
of $D$. This map has a right inverse, viz., the map
\begin{equation}
\label{def-varf}
F:H_0(\immpsr k n)\to H_0(\immsr k n)
\end{equation} that forgets the marked point.

In the following we will consider the combination of \eqref{hzero1}
with $\varphi^*$ thus obtaining a map:
\begin{equation}
i_k:H^{k(n-3)}(\imbr n)\to H^0(\immsr k n).
\label{hzero}
\end{equation}
A class $\omega$ in $H^0(\immpsr k n)$ will be called
equivariant if  $F^*\varphi^{*}\omega=\omega$.
\end{Rem}

Classes in $H^{k(n-3)}(\imbr n)$
can be constructed via trivalent graphs, as
shown in the sequel of this paper.
These classes have been firstly considered
in the 3-dimensional case, in connection with  perturbative
Chern--Simons quantum field theory.

\section{The Bott--Taubes construction}
\label{bott-taubes}

\subsection{Configuration spaces}
For any compact manifold $M$, we consider first the configuration
space $C_q^0(M)$ $\triangleq M^q\setminus\{\bigcup_S \Delta_S\}$,
where $S$ runs over the ordered subsets of the first $q$ integers
with $|S|\geq 2$, and $\Delta_S$ denotes  the (multi)-diagonal
labelled by $S$, namely, the subset of $M^q$ defined by the
equations $x_{j_1}=x_{j_2}=\cdots= x_{j_{|S|}}$, $j_i\in S$.

We consider then the compactification $C_q(M)$ of
$C_q^0(M)$ introduced in \cite{AS} as a modification of the
Fulton--MacPherson construction \cite{FMP}, as described below.

One has an obvious
inclusion of $C_q^0(M)\subset M^q$ and, for each diagonal
$\Delta_S$, one has a projection $C_q^0(M)\to
Bl(M^{|S|},\Delta_S)$ where $Bl$ denotes the differential-geometric
blowup (i.e., one replaces the given diagonal
$\Delta_S$ by
the sphere bundle of its normal bundle). This gives an imbedding
$C_q^0(M)\incul M^q\times \prod_{|S|\geq 2}Bl(M^{|S|},\Delta_S)$. The
space
$C^q(M)$ is then defined as the closure of  $C_q^0(M)$ in
the above space. The main fact about this compactification of
configuration spaces (see \cite{BT}) is the following:
\bth
The spaces $C_q(M)$ are smooth manifolds with corners, and all
the projections $C_q^0(M)\to C_{q-k}^0(M)$ extend to smooth projections
on the corresponding compactified spaces.
\label{thm-BT1}
\eth

The boundaries of $C_q(M)$ correspond
to the ``collision'' of at least
two of the $q$ points of $M$. Boundaries are the union of
different strata corresponding to the different ways in which all the
points may collide.
More precisely, let $S\subset\{1,\cdots,q\}$ be the labels of the colliding
points.
Let us insert in $S$ different levels of parentheses so that each pair of
parentheses contains at least two elements. Points in $M$ ``collide
at the same speed'' if they
belong to the same level of parentheses (points are assumed to
``collide'' starting
from the innermost parentheses).
The codimension of a given stratum is equal to the
number of pairs of parentheses.

We are mainly interested in codimension-1 strata, namely, in those strata
with no
internal parentheses. For these strata,
one calls {\em hidden faces}\/ those
corresponding to subsets $S$ with $|S|\geq 3$ and  {\em principal faces}\/
those for which $|S|=2$.

\subsubsection{The case of $S^1$}
If we choose $M$ to be $S^1$, then $C_q^0(S^1)$ is not connected.
We choose a connected component by fixing an order of the points on $S^1$
(consistent with its orientation).
It is then easy to see that this connected component is
given by $S^1\times\Sigma^0_{q-1}$ where $\Sigma^0_{q-1}$
is the ordinary open
$(q-1)$-dimensional simplex. We denote the closure of the connected component
of $C_q^0(S^1)$ by the symbol $C_q$.
This is given by the Cartesian product
of $S^1$ times a space $W_{q-1}$ that is obtained from the standard
closed $(q-1)$-simplex by a sequence of blowups (see the explicit
description in \cite{BT}).

\subsubsection{The case of $\reali^n$}
In the following we need a suitable compactification of $C_q^0(\reali^n)$.
Since $\reali^n$ is not compact, we cannot rely directly on the preceding
construction.

Instead, following \cite{BT}, we identify $\reali^n$ with
$S^n\setminus\{\infty\}$ and define $C_q(\reali^n)$ as the fiber
over $\infty\in S^n$ of $C_{q+1}(S^n)\to S^n$ (say,
projecting to the last factor).

This way, we also have a compactification (and corresponding boundary faces)
at infinity. (For example, $C_1(\reali^n)$ is the $n$-dimensional ball.)

\subsubsection{The case of an imbedding of $S^1$ into\/ $\reali^n$}
This is the case of interest for the rest of the paper.

Again following \cite{BT}, we define the
space $C_{q,t}(\reali^n)$ of $q+t$ distinct points in $\reali^n$, the first
$q$ of which are constrained on an imbedding of $S^1$,
as a pulled-back bundle as follows:
\begin{equation}
\begin{array}{ccc}
C_{q,t}(\reali^n) & \stackrel{\hat {ev}}{\lora} & C_{q+t}\,(\reali^n) \\
\downarrow \,p_1  &  & \downarrow \\
C_{q} \, \times \, \imbr n & \stackrel{ev}{\lora} & C_{q}\,(\reali^n)  \\
\end{array}
\label{puno}
\end{equation}
where the map $ev:C_{q} \, \times \, \imbr n \to C_{q}\,(\reali^n)$ is
the evaluation map applied to  $q$ distinct points in $S^1$ and
$\hat{ev}$ is its lift.

The diagram is commutative by construction. The main result
is the following theorem proved in \cite{BT}:
\bth
The spaces $C_{q,t}(\reali^n)$ are smooth manifolds with corners.
Moreover, the map $\Hat{ev}$ and the projection $p_1$ (and, more
generally, all projections $C_{q,t}(\reali^n)\to C_{q-k,t-l}(\reali^n)$)
are smooth.
\label{thm-BT2}
\eth

\subsection{Tautological forms}
It is not difficult to check that
the maps $\phi_{ij}: C_q^0(\reali^n) \to S^{n-1}$,
\[
\phi_{ij}(x_1,\ldots , x_q)\triangleq\frac{x_j-x_i}{|x_j-x_i|},
\]
extend to smooth maps on $C_q(\reali^n)$. In fact, it is enough
to consider the case $q=2$ and then apply Theorem~\ref{thm-BT1}.

Next we consider  the so-called {\em tautological forms}, which
are smooth as a consequence of Theorem~\ref{thm-BT2}. They are defined by
(see \cite{BT}):
\begin{equation}
\theta_{ij}(v^n)\triangleq \hat{ev}^*\phi_{ij}^*v^n\in
\Omega^{n-1}\left(C_{q,t}(\reali^n)\right),
\label{thetaij}
\end{equation}
where $v^n$ is a given normalized symmetric
smooth top form on $S^{n-1}$.

Other forms on $C_{q,t}(\reali^n)$ that we want
to consider are obtained by pulling back the  symmetric form $v^n$
via the map given by the combination of $p_1$ (s. \eqref{puno}) with
the map
\[
C_{q} \, \times \, \imbr n  \stackrel{pr_i \times id}{\lora}
S^1 \times \imb {\reali^n }.
\]
where $pr_i:C_q\to S^1$ denotes  the $i$th
projection.

The pullback of forms on $S^1\times \imbr n$ are forms on
$C_{q,t}(\reali^n)$.
The main example that we  have in mind is
the ``tangential tautological form''
\begin{equation}
\theta_{ii}(v^n)\triangleq (ev_i\circ D)^*v^n,
\label{thetaii}
\end{equation}
where $D$ is the normalized derivative
and  $ev_i = ev\circ(pr_i\times id)$.

\subsubsection{General properties of tautological forms}
\label{proptaut}
Taking into account  the definition
of the maps $\phi_ {ij}$ and of the tautological forms
(\ref{thetaij},~\ref{thetaii}), we have the following relations:
\begin{align}
\theta_{ij}(v^n)&=(-1)^ {n}\,\theta_{ji}(v^n),
\qquad i\not=j,\label{proptaut1}\\
\theta_{ij}(v^n)\,\theta_{uv}(v^n)&=(-1)^ {n+1} \,\theta_{uv}(v^n)\,
\theta_{ij}(v^n),
\label{proptaut2}\\
\theta_{ij}^2(v^n)&=0.   \label{proptaut3}
\end{align}
The first relation is due to the action of the antipodal map
on $S^{n-1}$,
the second relation is a consequence of the degree of the
tautological forms, and the third relation is an obvious consequence
of the fact that the square of a top form is zero.

Finally,
it may also be recalled that
the cohomology classes of the tautological forms generate the
whole cohomology of the configuration spaces of $\bbR^n$.

\subsection{Forms on the space of imbeddings}\label{ssec-fimb}
In order to have differential forms on $\imbr n$ we consider the
``push-forward,'' or fiber-integration. For any bundle ($p:E \to
B$) such that the fiber $F$ is an $m$-dimensional compact oriented
manifold (possibly with boundaries or corners), we define a map
$p_*$ from the space of $(p+m)$-forms on $E$ to the space of
$p$-forms on $B$, as follows:
\[
p_* \omega (X_1, \ldots ,X_m)\triangleq  \int_{F} \omega
(\tilde{X}_1, \ldots ,\tilde{X}_{m}, \cdot)
\]
where $\omega$ is a $(p+m)$-form on $E$ and $\tilde X_i$ is a
vector field on $B$ whose projection yields the vector field
$X_i$. The definition of $p_*$ is independent of the choice of the
lifts $\tilde X_i$.

{}From the sequence of maps:
\begin{equation}
\begin{array}{clc}
C_{q,t}(\reali^n) & &\\
\downarrow \,p_1& & \\
C_{q} \, \times \, \imbr n  & & \\
\downarrow \,p_2& & \\
\imb {\reali^n }& &
\end{array}
\nonumber
\end{equation}
we obtain, by fiber-integrating products of $\theta_{ij}(v^n)$s,
forms  on $\imbr n$ which are not necessarily closed since the
fiber is a manifold with corners. {}From the product of $k$
tautological forms  we obtain  a $((n-1)k-nt-q)$-form on $\imbr n
$.

\begin{Rem}
\label{rem-diff} Forms on $\imbr n$ obtained this way are
$\mathit{Diff}^+(S^1)$-equiv\-ariant. Observe in fact that an
orientation-preserving diffeomorphism of $S^1$ induces an
orientation-preserving diffeomorphism of $C_q$. Horizontality
follows then directly from the fiber integration along
configuration spaces of $S^1$, while invariance is a consequence
of the usual invariance of integrals under reparametrizations.
\end{Rem}

The exterior derivative of pushed-forward forms is given in terms
of the  generalized Stokes formula:
\[
d\; p_* \omega (X_1, \ldots ,X_m)=
p_* d \omega (X_1, \ldots ,X_m)+
(-1)^{\deg p_*^{\partial} \omega }
\,p_*^{\partial} \omega (X_1, \ldots ,X_m).
\]
The coboundary operator $d$ on the l.h.s.\ refers to the space
$\imb {\reali^n }$, while the coboundary operator on the r.h.s.\
refers to the space $ C_{q,t}(\reali^n)$. Moreover,
$p_{\ast}^{\partial} \omega$ is
given by
\[
p_*^{\partial} \omega(\gamma) \triangleq \int_{\partial
C_{q,t}(\reali^n,\gamma)} \omega,
\]
where $\partial C_{q,t}(\reali^n,\gamma)$
is the union of all the boundaries of codimension-1 of the fiber
over the imbedding $\gamma$.

If we denote by $\lambda$ a product of tautological forms,
then $d\lambda=0$.
So we have
\begin{equation}
d\, p_* (\lambda)=(-1)^{\deg p_*^{\partial} \lambda }\,
p_*^{\partial} (\lambda). \label{dp=pd+r}
\end{equation}
In Appendix~\ref{hidden} we will consider these boundary push-forwards
explicitly and show that, for $n>3$, only principal faces contribute.

\begin{Rem}
Let us consider the $j$th projection $p_j:C_q(\reali^n)\to
C_{q-1}(\reali^n)$, and let us define
\[
\tau_{ik} = p_{j*}\theta_{ij}(v^n)\,
\theta_{jk}(v^n)\in\Omega^{n-2}(C_{q-1}(\bbR^n)).
\]
As a consequence of \eqref{proptaut1}, \eqref{proptaut2} and
\eqref{dp=pd+r}, $\tau_{ik}$ is closed. But the $(n-2)$-nd cohomology
group of $C_q(\reali^n)$ is trivial. So the form $\tau_{ik}$ is exact.
\label{rem-triv}
\end{Rem}
\begin{Rem}
Another particular case is the integral over $C_q(\reali^n)$, $n>3$,
of a product of tautological forms with the condition that the situation
of the preceding Remark never occurs (that is, we assume that for each
point $i$, there are at least three tautological forms $\theta_{i\bullet}(v^n)$).
In this case, the result must be a number, and this will not vanish only
if the form degree matches the dimension of the space.

However, it is easy to prove that the form degree minus the dimension
of the configuration space is always greater or equal to $(n-3)q/2$.
So these integrals always vanish.
\label{rem-trivRn}
\end{Rem}

\section{The complex of decorated graphs in any dimension}
\label{ddco}
Push-forwards of products of tautological forms along configuration spaces
can be given a nice description in terms of graphs with a distinguished
oriented loop. In the following, we will always represent
the distinguished loop by a circle.

The idea is to represent each point in the configuration space as a vertex
of a graph with the convention that all vertices constrained on the
imbedding are put in order on the circle. Each tautological
form will then be represented by an edge not belonging to the circle.
(Actually, in the following we will reserve the term {\em edge}\/ only
to this kind of edges.)

In view of Remark~\ref{rem-triv},
we can restrict ourselves to considering only graphs whose vertices
not on the circle are at least trivalent.
Moreover, thanks to Remark~\ref{rem-trivRn} and to \eqref{proptaut3},
only connected graphs without multiple edges
may yield nonzero results.

To keep track of the orientation of the configuration space and of
the order in which one takes the product of tautological forms,
the vertices and the edges must be numbered. Moreover, to distinguish
between $\theta_{ij}(v^n)$ and $\theta_{ji}(v^n)$, one has to orient the edges.

However, thanks to the properties \eqref{proptaut1} and \eqref{proptaut2},
the decoration of graphs can be simplified, as will be explained
in subsection~\ref{decodiagr}.

The differential of a form on $\imb{\reali^n}$ will be related by
\eqref{dp=pd+r} to other push-forwards of products of tautological
forms. As explained in Appendix~\ref{hidden}, also these push-forwards
can be described in terms of graphs.
As a consequence, we may relate
the exterior derivative on the space of imbeddings to a certain coboundary
operator on the complex of graphs. This is explained in
subsection~\ref{ssec-cob}.

The whole construction will finally be summarized in
subsection~\ref{ssec-morph}, where we will also establish the
relation between the graph cohomology and the de~Rham cohomology
of the space of imbeddings.

\subsection{Decorated graphs in odd and even dimensions}
\label{decodiagr}
Following the above discussion,
we will consider {\em connected graphs}\/ consisting of an
{\em oriented circle}\/ and many {\em edges}\/ joining vertices
which may lie either on the circle ({\em external vertices}) or off
the circle ({\em internal vertices}). We also require that each
internal vertex should be at least trivalent.

In a graph we define a {\em small loop}\/ to be an edge whose
end-points are the same vertex. We call  a small loop external or
internal according to the nature of the corresponding vertex.
(External small loops will represent forms $\theta_{ii}(v^n)$ as
defined in \eqref{thetaii}, and internal small loops will be ruled
away by \eqref{equiv-isl}.)

Next we assign a {\em decoration}\/ to each graph
in order to take into account the specific properties
of the tautological forms:
\begin{itemize}
\item If $n$ is odd, then we label both
internal and external vertices
and assign  an orientation
(represented by an arrow) to each edge. We assume that the labelling
of the external vertices is cyclic w.r.t. the orientation of the circle.
Moreover, whenever we have an external small loop, we
fix an ordering of the two half-edges that form it;
notice that this ordering is chosen independently from the edge orientation.
\item If $n$ is even, then the
decoration consists in the labelling of the external vertices
and of all the edges. Again we assume that the labelling
of the external vertices is cyclic w.r.t. the orientation of the circle.
\end{itemize}

We now define $\mathcal{D}^\prime_o$ ($\mathcal{D}^\prime_e$)
to be the real vector space generated by decorated graphs of odd (even)
type (some examples of elements of these spaces are in
figures~\ref{figura2}, \ref{figura5}, \ref{figura3} and~\ref{figura4}).

As explained at the beginning of the Section, we actually do not need
the whole spaces of graphs. We will restrict ourselves to the interesting
spaces by dividing $\mathcal{D}^\prime_o$ and $\mathcal{D}^\prime_e$
by certain equivalence relations.

The first two relations do not depend on the decoration and are as follows:
\begin{align}
\Gamma&\sim 0,\ & &\text{if
 two vertices in $\Gamma$ are joined by more than one edge,}
\label{equiv-de}\\
\Gamma&\sim 0,\ & &\text{if $\Gamma$ contains an internal small loop.}
\label{equiv-isl}
\end{align}
(The first relation is motivated by \eqref{proptaut3}, and the second
by the fact that we cannot associate to an
internal small loop any tautological form.)

Next, for any given pair of graphs $\widehat\Gamma$ and $\Gamma$
that differ
only for the decoration, we introduce the following equivalence
relations:
\begin{itemize}
\item For $\Gamma,\widehat\Gamma\in\mathcal{D}^\prime_o$,
\begin{equation}
\Gamma\sim (-1)^ {\pi_1+\pi_2+l+s}\, \, \widehat\Gamma,
\label{relazd}
\end{equation}
where $\pi_1$ is the order of the permutation of the internal vertices,
$\pi_2$ is the order of the (cyclic) permutation of external vertices,
$l$ is the number of edges whose orientation has been reversed,
and $s$ is the number of external small loops on which the ordering
of the half-edges has been reversed.

\item For $\vrule width 0pt height 20pt
\Gamma,\widehat\Gamma\in\mathcal{D}^\prime_e$,
\begin{equation}
\Gamma\sim(-1)^ {\pi+l}\,\,\widehat\Gamma,
\label{relazp}
\end{equation}
where $\pi$ is the order of the (cyclic) permutation of the external
vertices, and  $l$ is the  order of the permutation of the edges.
\end{itemize}
In order to have a well-defined, one-to-one correspondence between decorated
graphs and the
push-forwards of tautological forms as described at the
beginning of  the Section,
we need to quotient
$\mathcal{D}^\prime_o$ and
$\mathcal{D}^\prime_e$ with respect to
the equivalence relations (\ref{equiv-de},\ref{equiv-isl},\ref{relazd})
and, respectively, by (\ref{equiv-de},\ref{equiv-isl},\ref{relazp}).
Namely, we define:
\[
\mathcal{D}_o:=\mathcal{D}^\prime_o/\sim
\quad\text{and}\quad
\mathcal{D}_e:=\mathcal{D}^\prime_e/\sim.\]

\subsubsection{Order and degree of decorated graphs}
The {\em order}\/ of a graph $\Lambda$ (i.e., minus its Euler
characteristic) is defined as
\begin{equation}
\ord\Lambda = e - v_i,
\label{deford}
\end{equation}
where $e$ is the number of edges and $v_i$ is the number of internal
vertices.

The {\em degree}\/ of a graph $\Lambda$ is defined as
\begin{equation}
\deg\Lambda = 2 \, e - 3 \, v_i - v_e,
\label{defdeg}
\end{equation}
where $e$ is the number of edges, $v_e$ is the number of
external vertices
and $v_i$ is the number of internal vertices.

In the particular case when
the graph has only trivalent internal vertices and
univalent external vertices, its degree is zero and its order is half
the total number of vertices.

We consider $\mathcal{D}_o$ and
$\mathcal{D}_e$ as
{\em graded vector spaces}\/ with respect to both the {\em order}\/ and the
{\em degree}.

We denote by $\mathcal{D}^{k,m}_o$
and $\mathcal{D}^{k,m}_e$ the
equivalence classes of decorated
graphs of {\em order $k$ and degree $m$.}

\subsection{A coboundary operator for decorated graphs}\label{ssec-cob}
Now we want to introduce a coboundary operator on each space
$\mathcal{D}_o$ and $\mathcal{D}_e$.

As explained at the beginning of this Section, we actually look for
coboundary operators that, under the correspondence between graphs
and configuration space integrals, are related to the exterior derivative on
$\imbr n$.

We will first define these operators on $\mathcal{D}^\prime_o$
and $\mathcal{D}^\prime_e$, and then prove that they descend to the
quotients. These operators (both on the primed space and on their quotients)
will be denoted by $\delta_o$ and $\delta_e$ respectively.
When considering graphs of unspecified parity,
we will simply use the symbol $\delta$.

First of all we introduce some terminology.
\bdefiniz
We call {\em chord}\/ an edge whose end-points are distinct
external vertices
and {\em short chord}\/ a chord whose end-points are consecutive
vertices on the circle.
We call {\em regular edge}\/ an edge that is neither a chord nor a
small loop.
Finally we call  {\em arc}\/ a
portion of the circle bounded by two consecutive external vertices.
\edefiniz

For any graph $\Gamma$,
$\delta\Gamma$ will be, by definition, a signed sum
of decorated graphs obtained by contracting, one at a time, all
the regular edges and all the arcs of $\Gamma$.
Notice that the contraction of
an arc joining the vertices of a short
chord will produce an external small loop.
In the odd case, we order its
half-edges consistently with the orientation of the circle.

Edges and vertices are then relabelled as follows after
contraction: if the new graph is obtained by contracting vertex
$i$ with vertex $j$, then we relabel the vertices by lowering by
one the labels of the vertices greater than max($i,j$) and assign
the label min($i,j$) to the vertex where the contraction has
happened. If we contract the edge $\alpha$, we lower  by one the
labels of all the edges greater than $\alpha$.

Moreover, we associate to each
contraction  a sign defined as follows:
\begin{itemize}
\item On the space  $\mathcal{D}_o'$,
the sign associated to the contraction of the edge or of the arc
joining the vertex $i$ to the vertex $j$
(using the orientation of the edge or the arc)
is given by
\begin{equation}
\sigma(i,j)=\left\{
\begin{array}{lcl}
(-1)^ {j}& &\mbox{if}\,\, j>i,\\
 (-1)^ {i+1}& &\mbox{if}\,\,j<i.
\end{array}
\right.
\label{sigma1}
\end{equation}
\item On the space $\mathcal{D}_e'$,
the sign associated to the contraction of the arc joining vertex $i$
to vertex $j$ is given by
\begin{equation}
\sigma(i,j)=\left\{
\begin{array}{lcl}
(-1)^ {j}& & \mbox{if}\,\, j>i,\\
 (-1)^ {i+1}& & \mbox{if}\,\,j<i;
\end{array}
\right.
\label{sigma2}
\end{equation}
while, if we contract edge $\alpha$, we have the following sign:
\begin{equation}
\sigma(\alpha)=(-1)^{\alpha+1+v_e},
\label{sigmaedge}
\end{equation}
where $v_e$ is the number of external vertices.
\end{itemize}

The following Theorem shows that for both odd and even case,
$\mathcal D$ is a complex of graded vector spaces.

\bth
The operators $\delta_o$ and $\delta_e$
descend to $\mathcal{D}_o$ and, respectively, $\mathcal{D}_e$.
They are both coboundary operators there, and we have
\[
\delta_o: \mathcal{D}^{k,m}_o \to \mathcal{D}^{k,m+1}_o\,\,\,\quad
\delta_e: \mathcal{D}^{k,m}_e \to \mathcal{D}^{k,m+1}_e.
\]
\eth
The corresponding cohomology groups will be denoted in the following by
$H^{k,m}(\mathcal{D}_o)$ and $H^{k,m}(\mathcal{D}_e)$.

\begin{proof}
First we consider graphs of odd type and review the proof given in
\cite{BC1}.

Let us consider the contraction of $(ij)$,
i.e., of the edge or portion of circle between $i$ and
$j$. If we exchange
$i$ and $j$ or reverse an arrow, we get a minus sign;
in both cases the roles of $i$ and $j$ are interchanged and
we have $\sigma(i,j)
=-\sigma(j,i)$. Therefore $\delta_o$ is compatible with such an exchange.

Let us choose another vertex $k$, and exchange $j$ and $k$. We can
assume $i<j$ and that $(ij)$ is oriented from $i$ to $j$. First we
suppose $k>i$ and $k>j$. If we contract $(ij)$ we get a factor
$(-1)^j$; if we exchange $j$ and $k$ and then contract $(ik)$ we
get a factor $(-1)^{k+1}$. Obviously the underlying graph is the
same in the two cases. We want to prove that the relabelling of
one the two decorated graphs yields the other one, i.e., the two
decorated graphs define the same element in $\mathcal {D}_o$. The
indices lowered by one are, in the first case, all those greater
than $j$ and in the second case all those greater than $k$. The
set of vertices that, in the first case, are labelled as
$j,(j+1),\cdots,(k-2),(k-1)$, are labelled, in the second case, as
$(j+1),(j+2),\cdots, (k-1),j$ and the sign of the relevant
permutation is  $(-1)^{k-j+1}$.

In summary:
\begin{itemize}
\item if we contract $(ij)$ we get a sign $(-1)^j$;
\item if we swap $k$ and
$j$, contract $(ik)$ and then relabel to get the same graph as in the
previous
case, we get the sign $(-1)^{k+1}(-1)^{k-j+1}=(-1)^{j}$.
\end{itemize}

Similarly, we can treat the case $k<i$. All other cases (contraction of
$(ij)$ and swapping of $l$ with $k$) are trivial.

Now let us prove that $\delta_o^2=0$
by showing that contracting two different
pairs $(ij)$ and $(rs)$ in opposite order yields the opposite sign.

If we have $i\neq r$ and $j\neq s$, then we can always assume that
$i<j$, $r<s$ and $j<s$. Contracting $(ij)$ gives a sign $(-1)^j$ and lowers $s$
by one, so contracting $(rs)$ gives $(-1)^{s-1}$. If we contract
$(rs)$ first, then $j$ is {\em not} lowered by one and the global sign is
$(-1)^s(-1)^j$.

If $s=j$ and $i\neq r$, we can pass from a (double)-contraction to the
other one by exchanging $i$ and $r$, with a change in  sign in
$\mathcal{D}_o$. The same holds for $s\neq j$ and $i=r$.\\

Next let us consider $\mathcal{D}^\prime_e$.

Again we have to show that a permutation of external vertices or
of edges does not affect $\delta_e$. The case when we swap
external vertices is identical to the case of $\mathcal{D}_o$, so
we just have to verify what happens when we swap two edges
labelled by $\alpha$ and $\beta$, with $\alpha<\beta$.

We claim that if we contract the edge $\alpha$
we obtain the same result as if we swap the edge $\alpha$ with
the edge $\beta$ and subsequently contract the edge $\beta$.
In fact in the first case the result  is $(-1)^{\alpha+1+v_e}$ 
times a graph whose edges are:
$(1,\ldots, \alpha-1,\alpha+1,\ldots,\beta-1,\beta,\beta+1,\ldots,t)$,
while in the second case the result is $(-1)^{\beta+v_e}$ times a graph
whose edges are $(1,\ldots, \alpha-1,\alpha+1,\ldots,\beta-1,\alpha,
\beta+1,\ldots,t)$. We now permute the labels of the edges of
this last graph and obtain  $(1,\ldots, \alpha-1,\alpha,\alpha+1,
\ldots,\beta-1,\beta+1,\ldots,t)$; this
permutation has order $(\beta-1)-(\alpha+1)+1=\beta-\alpha-1$.
The total sign is therefore:
$(-1)^{\beta+v_e}(-1)^{\beta-\alpha-1}=(-1)^{\alpha+1+v_e}$, i.e.,
the same result obtained by contracting the edge $\alpha$.

Finally, we have to show that $\delta_e^2=0$. As in the odd case,
the proof consists in showing that if we make two contractions
in different order, then we have opposite signs and hence
the relevant graphs cancel.

This is obviously true if we contract two pairs of external
vertices or two edges. Now we contract the arc  between
two consecutive external vertices (say $i$ and $j$, with $i<j$)
and an edge $\alpha$. Remember that the number $v_e$ of external
vertices appears in equation~\eqref{sigmaedge}. If we contract
$\alpha$ first, we do not change the labels of the external
vertices and we get the global sign $(-1)^{\alpha+1+v_e}(-1)^j=
(-1)^{\alpha+1+v_e+j}$. If we contract $(ij)$ first, then the
number of edges is lowered by one and so the global sign is
$(-1)^j(-1)^{\alpha+1+(v_e-1)}=(-1)^{\alpha+v_e+j}$.
\end{proof}

\subsection{The configuration space integral as a morphism of complexes}
\label{ssec-morph}
We can finally give all the details of the construction described
at the opening of the Section.

First of all, we fix $n>3$ and choose a form $v^n$ through which
we define our tautological forms.

\begin{Def}
\label{def-gf}
Let $\alpha$ be the antipodal map on $S^{n-1}$. We call {\em
symmetric form} a normalized element $v^{n}\in\Omega^n(S^{n-1})$
that satisfies $\alpha^*v^{n} = (-1)^{n}\,v^{n}$.
\end{Def}
Obviously the standard volume form on $S^{n-1}$ is a symmetric
form and, moreover, no nontrivial top form $w^n$ on $S^{n-1}$ can
satisfy the condition $\alpha^*w^{n} = (-1)^{n-1}\,w^{n}$.

From now on we will always assume that $v^n$ is a fixed normalized
symmetric form on $S^{n-1}$.

Then we associate to every class of graphs $\Gamma$ a form $I(v^n)(\Gamma)$
on $\imbr n$
in terms of a configuration space integral as follows:
\begin{itemize}
\item In the odd case, the edge 
joining the vertex $i$ to the vertex $j$
is replaced by the form $\theta_{ij}(v^n)$. The external
small loop at $k$ is replaced by $(-1)^{l+s}\,\theta_{kk}(v^n)$, where
$l=0$ if the edge orientation
agrees with the ordering of the half-edges and $l=1$ otherwise, while
$s=0$ if the ordering of the half-edges corresponds to
the orientation of the circle and $s=1$ otherwise.
Since all the forms are even, we do not have to say in which order
we take their product. The orientation of the configuration space
is determined by the numbering of the vertices.
\item In the even case,
we first have to choose a numbering of the internal vertices as well.
Then the edge joining the vertex $i$ to the
vertex $j$ is replaced by the form $\theta_{ij}(v^n)$, and an external small
loop at $k$ is replaced by the form $\theta_{kk}(v^n)$.
The numbering of the edges tells us in which order we have to take the
product of tautological forms,
while the numbering of the external vertices determines
the orientation of the configuration space. (The numbering of internal
vertices is on the other hand irrelevant.)
\end{itemize}
We will denote by $I(v^n)$ the linear extension to $\mathcal D$ of the
map 
just described.
Then we have the following:
\bth
For any $n>3$ and for any symmetric form $v^n$,
\begin{align*}
I(v^n)&\colon(\mathcal D_o^{k,m},\delta_o)
\lora (\Omega^{m+(n-3)k}(\imbr n),d),\,
&&\text{ if $n$ is odd,} \\
I(v^n)&\colon(\mathcal D_e^{k,m},\delta_e)
\lora(\Omega^{m+(n-3)k}(\imbr n),d),\,
&&\text{ if $n$ is even,}
\end{align*}
are chain maps.\label{thm-Ivn}
\eth
\begin{proof}
The coboundary operator $\delta$ has been defined in subsection~\ref{ssec-cob}
in such a way that it corresponds to the coboundary operator $d$ via Stokes'
Theorem if one considers only principal faces
(see Appendix~\ref{hidden} and in particular
Thm.~\ref{thm-princ} and Rem.~\ref{rem-princ}).
The fact that $I$ is actually a chain map then follows from Theorem
\ref{thm-vt} according to which we can neglect hidden faces.

As for the degree of these maps, this is easily
computed: if the graph $\Lambda$ has $e$ edges,
$v_i$ internal vertices and $v_e$ external vertices, then, by \eqref{defdeg},
$\deg \Lambda= 2 \, e - 3 \, v_i - v_e $.
On the other hand the degree of the corresponding
differential form is
\[
\deg I(\Lambda)= (n-1) \, e - n \, v_i - v_e =
\deg \Lambda +(n-3)\ord \Lambda,
\]
where $\ord \Lambda$ is the order of the graph as defined in
\eqref{deford}.
\end{proof}
In the following we will denote also by $I(v_n)$ the map induced in cohomology:
\begin{align}
\label{oddmorph}
I(v_n)&\colon H^{k,m}(\mathcal D_o)
\lora H^{m+(n-3)k}(\imbr n),\,
&&\text{ if $n$ is odd,} \\
I(v_n)&\colon H^{k,m}(\mathcal D_e)
\lora H^{m+(n-3)k}(\imbr n),\,
&&\text{ if $n$ is even.}
\label{evenmorph}
\end{align}
The case $m=0$ is particularly interesting and will be discussed in the next
section. In Section~\ref{nontriv} we will prove that in this case the above
homomorphisms are actually injective.

Recall finally that, as observed in Remark~\ref{rem-diff}, the image of
$I(v^n)$ lies in the subspace of
$\mathit{Diff}^+(S^1)$-equivariant forms. Thus, we can
produce elements of the $\mathit{Diff}^+(S^1)$-equivariant cohomology
of $\imbr n$.

\subsubsection{The dependency on the symmetric form $v^n$.}
We want now to consider the dependency of the chain map $I$ on the
choice of the symmetric form. We have the following
\bprop
Let $n>4$.
If $v_0^n$ and $v_1^n$ are two symmetric forms and\/ $\Gamma$ is a cocycle,
then $I(v_1^n)(\Gamma)-I(v_0^n)(\Gamma)$ is an exact form.
\label{prop-Iexact}
\eprop
\begin{proof}
Let us write $v_1^n-v_0^n=d w^n$, where $w^n\in\Omega^{n-2}(S^{n-1})$.
We can assume that $\alpha_n^*w^n=(-1)^n\,w^n$.
Then $v_t^n=v_0^n+t\,d w^n$, $t\in[0,1]$, interpolates between $v_0$ and $v_1$.
We now define
\[
\Tilde v^n \triangleq v_0^n + d(t\,w^n) \in \Omega^{n-1}(S^{n-1}\times[0,1]).
\]
This form is still closed and symmetric: $(\alpha_n\times id)^*\Tilde v^n=
(-1)^n\,\Tilde v^n$.

Denoting by $i_t:S^{n-1}\incul S^{n-1}\times[0,1]$ the inclusion at
$t\in[0,1]$, we also have
\begin{align*}
i_t^*\Tilde v^n &= v^n_t,\\
\int_{S^{n-1}} i_t^*\Tilde v^n &= 1.
\end{align*}

Using this extended symmetric form, we can define extended tautological forms
by
\[
\Tilde\theta_{ij}\triangleq (\hat{ev}\times id)^*(\phi_{ij}\times id)^*
\Tilde v^n\in
\Omega^{n-1}\left(C_{q,t}(\reali^n)\times[0,1]\right).
\]
If we now replace the edges of a graph by these extended tautological
forms, after integrating over the configuration space we will get a form
on $\imbr n\times[0,1]$. Denote by $\Tilde I(\Tilde v^n)$ this map.
Observe that, denoting by $j_t:\imbr n\incul
\imbr n\times[0,1]$ the inclusion at $t$, we have
\[
j_t^*\Tilde I(\Tilde v^n)(\Gamma) = I(v^n_t)(\Gamma).
\]
If $\Gamma$ is a cocycle, then the results of Appendix~\ref{hidden}
(in particular Thms.~\ref{thm-princ} and~\ref{vt01})
also
show that $\Tilde I(\Tilde v^n)(\Gamma)$ is a closed form in
$\Omega^*(\imbr n\times[0,1])$. As a consequence,
\begin{equation}
I(v^n_1)(\Gamma)-I(v^n_0)(\Gamma)= d\pi_*\Tilde I(\Tilde v^n)(\Gamma),
\label{Iv10n}
\end{equation}
where $\pi$ is the projection $\imbr n\times[0,1]\to\imbr n$, and we have used
again the generalized Stokes formula.
\end{proof}
Thus, for $n>4$ the homomorphisms \eqref{oddmorph} and \eqref{evenmorph} do
not depend on the choice of the symmetric form.

If $n=4$, then the homomorphisms might depend on the chosen symmetric form
$v^4$ (see Rem.~\ref{v4}).

Anyway, we will prove (see Thm.~\ref{HD-HImm})
that the integrals of
$I(v^n)(\Gamma)$, $\Gamma\in H^{*,0}(\mathcal D)$, on the cycles
given in \eqref{zyklon} do not depend on $v^n$, even for $n=4$.

\bconv
In the rest of this paper,
in all
the tautological forms, we will omit
the explicit dependency on the symmetric form $v^n$.
\econv

\section{Cocycles of trivalent graphs}
\label{sect-triv}
Trivalent graphs are defined as (decorated)
graphs having exactly
one edge for each external vertex, while exactly three edges merge into each
internal vertex. Notice that in particular trivalent graphs do not have
external small loops.
Moreover, trivalent graphs span $\mathcal{D}^{k,0}$.
 Particular cases of trivalent graphs are of course chord diagrams.

We will
look for cocycles that
are linear combinations of trivalent graphs, that is, elements of
$H^{k,0}(\mathcal{D})$.
We have the following
\bprop
If\/ $\Gamma=\sum_i c_i\Gamma_i\in H^{k,0}(\mathcal{D})$
is a cocycle given by a linear combination of
trivalent graphs\/ $\Gamma_i$, then at least one\/
$\Gamma_i$ is a chord diagram.
Moreover, no graph in\/ $\Gamma$  may contain a short chord.
\label{thm-trichord}
\eprop
\begin{proof}
Let $l$ be the  minimum number of internal
vertices among the graphs $\Gamma_i$.
The first statement is equivalent to saying that $l=0$.

In fact, assume on the contrary that $l>0$, and let $\Gamma_j$ be a
graph (which does not vanish in $\mathcal{D}$)
with exactly $l$ internal vertices. Since we consider only connected
graphs, there will be at least one internal vertex connected by one edge
(call it $f$) to an external vertex.

In $\delta\Gamma_j$ there will then be
a graph $\Gamma_j'$
obtained by contracting the edge $f$.

First of all, notice that $\Gamma_j'$ does  not vanish in $\mathcal{D}$.
In fact, if it did, then there would be an automorphism $\varphi$
of the graph
underlying $\Gamma_j'$ that would yield $-\Gamma_j'\in\mathcal{D}$.
Observe that the only 4-valent vertex in $\Gamma_j'$ has to be mapped to
itself by $\varphi$.
We can now extend this to an automorphism of $\Gamma_j$ by decollapsing $f$
after the application of $\varphi$ and
deciding that
$f$ is mapped into itself
and each of its end-points are mapped into themselves (notice that
we cannot interchange the end-points of $f$ since one is internal and the
other is external). So by the extended automorphism we would prove that
also $\Gamma_j$ vanishes in $\mathcal{D}$.

Since $\Gamma$ is a cocycle,
there must be other graphs such that their images under the application
of $\delta$ contain $\Gamma_j'$. But there are only
two possible graphs with
this property, namely, those obtained by splitting the unique four-valent
external vertex in $\Gamma_j'$ in the two possible ways. But both
these graphs have $l-1$ internal vertices, so $l$ cannot be the minimum.

To prove the second statement, observe first that $H^{1,0}(\mathcal D)=0$
(in the odd case since the chord diagram is not closed, in the even case
since the chord diagram vanishes by symmetry).

So assume $k>1$. This means that the number of vertices is greater than 2.
Since we consider only connected graphs, this means that, if
a graph $\Gamma_i$ contains a short chord, then
there is only
one arc (call it $a$) that has the same end-points as the short chord.

Let $\Gamma_i'$ be the graph in $\delta\Gamma_i$ obtained by contracting $a$.
Notice that $\Gamma_i'$ contains an external small loop.
Thus, since we do not allow internal small loops,
there is no other graph whose image under the application of $\delta$
may contain $\Gamma_i'$.

But as above we can prove that $\Gamma_i'$ does not vanish
(unless $\Gamma$ is zero itself.)

So $\Gamma$ cannot be a cocycle.
\end{proof}

Applying the homomorphisms \eqref{oddmorph} and \eqref{evenmorph}
to a trivalent cocycle $\Gamma$, we will then get cohomology classes
on $\imbr n$ of degree $(n-3)\ord\Gamma$.

The question whether nontrivial
cocycles in the graph complex represent  nontrivial
cocycles in the cohomology ring $H^*(\imbr n)$ will be addressed to
in the next Section.

In the rest of this Section, we will discuss some examples of trivalent
cocycles.

\subsection{The odd case}\label{ss-odd}
Trivalent graphs have been widely studied
(\cite{BN95,GMM,BT,Bott,AF}) in the case $n=3$. Here
all graphs are associated to
zero-forms, i.e., functions on $\imb{\reali^3}$) and,
at each order, it is possible to find $\delta$-closed
linear combinations of graphs.

Elements of $H^0(\imb{\reali^3})$
are nothing but  constant functions on connected components of
the space $\imbr 3$, i.e., knot invariants. Such invariants
are related to topological quantum field theories, namely,
to the Chern--Simons theory \cite{Wit} and to
$BF$ theories (see \cite{CCRFM95} and references therein).
Moreover, they are Vassiliev invariants (i.e., invariants of finite type)
(see \cite{V1,BN95,AF}).

When $n>3$, the homomorphism \eqref{oddmorph} implies that
we can construct cohomology classes on $\imbr n$
using exactly
the same linear combinations of graphs that give knot invariants
in three dimensions.

\begin{figure}[ht!]
\cl{\resizebox{8cm}{!}{\includegraphics{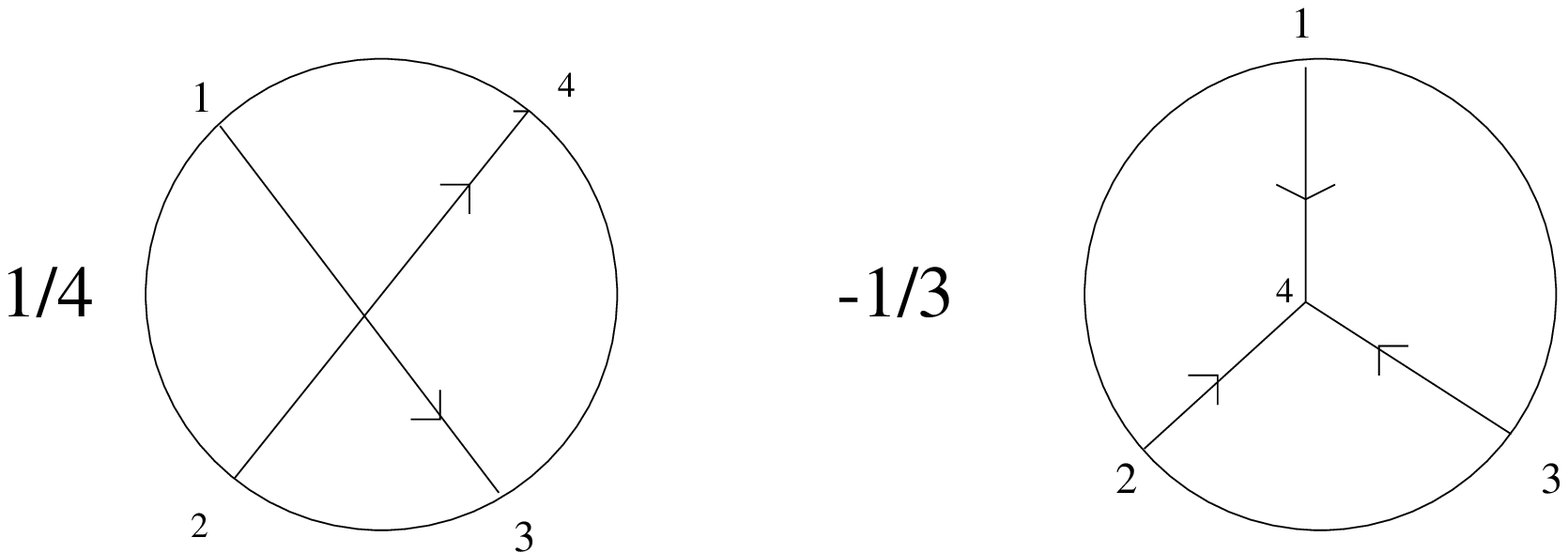}}}
\caption{Cocycle of odd type at order 2}
\label{figura2}
\end{figure}

For instance the simplest cocycle
\[
1/4\, \,\int_{C_{4}}\theta_{13} \theta_{24} -
1/3 \,\,\int_{C_{3,1}}\theta_{14} \theta_{24} \theta_{34},
\]
see figure \ref{figura2}, represents an element of
$H^{2n-6}(\imbr n)$.

\begin{figure} [ht!]
\cl{\resizebox{10cm}{!}{\includegraphics{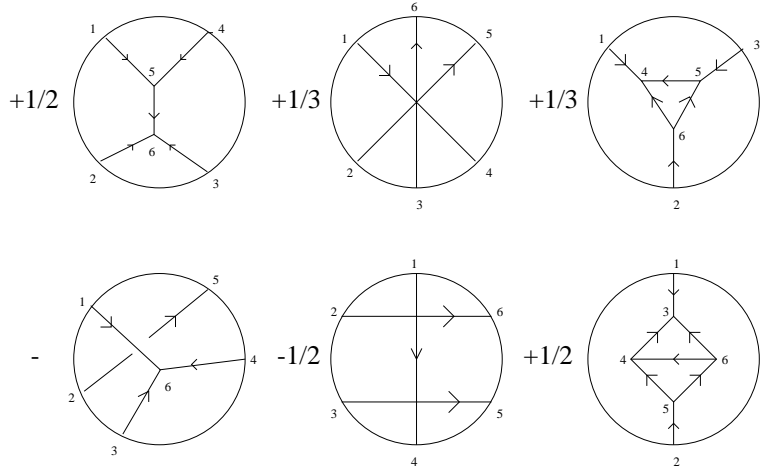}}}
\caption{Cocycle of odd type at order 3}
\label{figura5}
\end{figure}

At order 3, there is only one cocycle, which was calculated in
\cite{AL} and \cite{HS} (s. also \cite{BN95}). We
show it in figure~\ref{figura5}. In terms of integrals
of tautological forms, it is
given by:
\begin{multline}  \nonumber
\frac1{2}\,\int_{C_{4,2}} \theta_{15} \theta_{26} \theta_{36} \theta_{45}
\theta_{56}
+\frac13\,\int_{C_{6}} \theta_{14} \theta_{25} \theta_{36}
+ \frac13\,\int_{C_{3,3}} \theta_{14} \theta_{26} \theta_{35} \theta_{64}
\theta_{65} \theta_{54} +\\
-\, \int_{C_{5,1}}\theta_{16} \theta_{25} \theta_{36} \theta_{46}
-\frac12 \,\int_{C_{6}} \theta_{14} \theta_{26} \theta_{35}
+\frac12 \,\int_{C_{2,4}} \theta_{13} \theta_{25} \theta_{54} \theta_{56}
\theta_{64} \theta_{63} \theta_{43} .
\end{multline}
Recall now that every Vassiliev invariant in three dimensions
produces a nontrivial cocycle of trivalent graphs of odd type
\cite{BN95}. It is also well known that there are nontrivial
Vassiliev invariants at any order (e.g., coefficients of the
Alexander--Conway or of the Jones polynomials). Hence,
$H^{k,0}(\mathcal D_o)$ contains nontrivial elements for any
$k\geq 2$.

\subsection{The even case}\label{ss-even}
For $n$ even, we  use the homomorphism \eqref{evenmorph} to
construct cohomology classes of $\imbr n$.

\begin{figure}[ht!]
\cl{\resizebox{8cm}{!}{\includegraphics{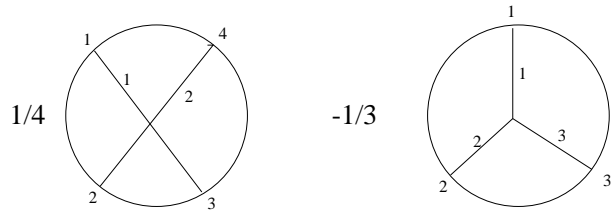}}}
\caption{Cocycle of even type at order 2}
\label{figura3}
\end{figure}

\begin{figure}[ht!]
\cl{\resizebox{10cm}{!}{\includegraphics{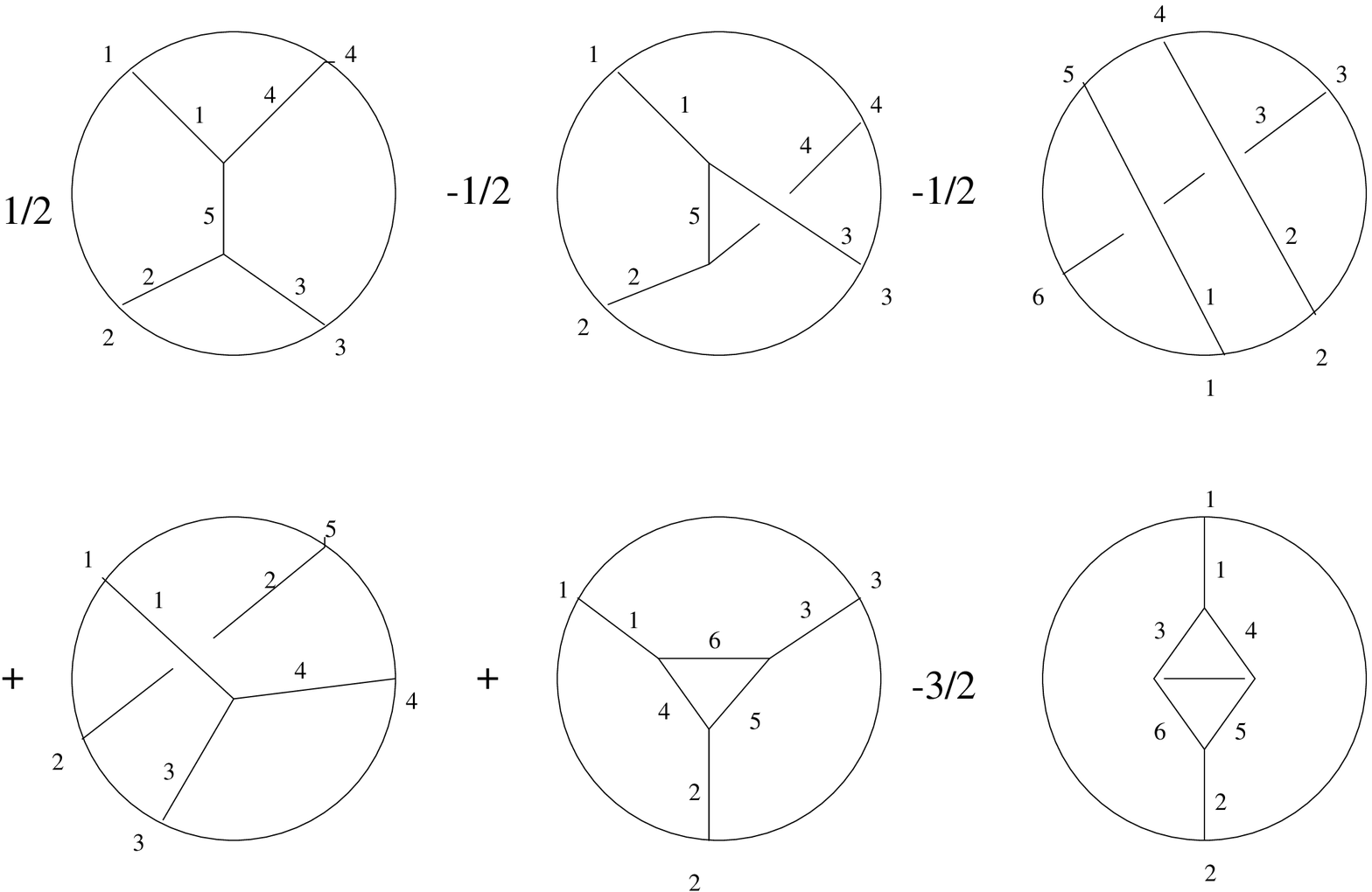}}}
\caption{Cocycle of even type at order 3}
\label{figura4}
\end{figure}

It is easy to show that, at order two, there is only
one cocycle in $\mathcal D_e$ (see figure~\ref{figura3}),
which induces, for every even $n$, the element
\[
1/4\,\, \int_{C_4}\theta_{13}\theta_{24} -
1/3 \,\,\int_{C_{3,1}}\theta_{14}\theta_{24} \theta_{34}
\in H^{2n-6}(\imb{\reali^n}).
\]
At order 3 the vector space generated by trivalent graphs of even type
has again dimension one, and the generator is given explicitly
in figure~\ref{figura4}.

The corresponding element of $H^{3n-9}(\imbr n)$ can be written as:
\begin{multline}   \nonumber
\frac12\,\int_{C_{4,2}} \theta_{1a} \theta_{2b} \theta_{3b} \theta_{4a}
\theta_{ab} -
\frac12\, \int_{C_{4,2}}\theta_{1a} \theta_{2b} \theta_{3a} \theta_{4b}
\theta_{ab}
-\frac12\,\int_{C_{6}} \theta_{15} \theta_{24} \theta_{36} + \\
+ \int_{C_{5,1}}\theta_{1a} \theta_{25} \theta_{3a} \theta_{4a} +
\int_{C_{3,3}}\theta_{1a} \theta_{2b} \theta_{3c} \theta_{ab}
\theta_{bc} \theta_{ca} -
\frac32\, \int_{C_{2,4}}\theta_{1a} \theta_{2b} \theta_{ac}
\theta_{ad} \theta_{bd}  \theta_{bc} \theta_{cd}.
\end{multline}
Here we have labelled with numbers the external vertices and with
letters the internal ones.

Moreover, we prove in \cite{CCRL2} that, as in the odd case, there
exist nontrivial cohomology classes of decorated trivalent graphs of
even type in an infinite number of orders $k$.

\section{Trivalent graphs and nontrivial classes in\nl
$H^*(\imbr n)$}\label{nontriv}

Let us now consider any nontrivial class in the cohomology
of graphs. Under \eqref{oddmorph} or \eqref{evenmorph},
this class will
represent a possibly trivial cohomology class
on $H^*(\imbr n)$.

This Section is devoted to prove that, however,
classes in the cohomology of trivalent graphs---that is, in
$H^{*,0}(\mathcal{D})$---yield
nontrivial classes in the cohomology of $\imbr n$.

The proof is based on the analysis of
Section~\ref{vorder} and particularly on the criterion provided by
Corollary~\ref{crithzero}.

The main result of the present Section is the following:
\bth
Let $n>3$. An $[(n-3)k]$th cohomology class on $\imbr n$, corresponding
via \eqref{oddmorph} or \eqref{evenmorph}
to a linear combination  of trivalent graphs, yields via
\eqref{hzero1} a nontrivial\/ $0$th cohomology class  on $\immpsr k n$.
\label{thm-nontriv}
\eth

\begin{proof}
Let $\Gamma=\sum_i\lambda_i\,\Gamma_i\in H^{k,0}(\mathcal{D})$
be the given cocycle of trivalent graphs. Let us denote
by $\omega$ ($\omega_i$) the image of $\Gamma$ ($\Gamma_i$) under
\eqref{oddmorph} or \eqref{evenmorph}. So
$\omega=\sum\lambda_i\omega_i$ is a closed $k(n-3)$-form  on $\imbr n$.

The forms $\omega_i$ are given by integrals of products of
tautological forms over the corresponding
configuration spaces $C_{q_i,t_i}(\reali^n)$ with
$q_i\leq 2k$, and
$q_i=2k$ if and only if $\Gamma_i$ is one of the chord diagrams
which are necessarily contained in $\Gamma$ by Prop.~\ref{thm-trichord}.

We integrate first over the {\em internal vertices} and denote
the result by $\mu_i$; then we have $\omega_i=
\int_{C_{q_i}} \mu_i$.

We consider $\gamma\in\immpsr k n$, and to each double point
$p^j,\;j=1,\dots, k$, of $\gamma$ we associate a ball $D_j$ of radius
$\epsilon>0$. The intersection $\gamma\cap D_j$ is given by
$l^j_1(\epsilon)\cup l^j_2(\epsilon)$, where $l^j_{a_j}(\epsilon)$,
$a_j=1,2$, are two closed segments as in
figure~\ref{figura1}.

We now define $C_{q_i}^{(\epsilon;j,a_j)}$
as the open subset
of $C_{q_i}$  such that the image through $\gamma$ of all
its projections over $S^1$ does not intersect  $l^j_{a_j}(\epsilon)$.
We also define
\[
C^{\epsilon}_{q_i}\triangleq \complement\left(
\bigcup_{j=1}^{k}\bigcup_{a_j=1}^2\; C_{q_i}^{(\epsilon;j,a_j)}
\right)
\]
In other words, this complementary set
$C_{q_i}^{\epsilon}$ is equal to the subset of $C_{q_i}$ for which
all the projections over $S^1$ yield, through $\gamma$,
one and only one element in each
$l^j_{a_j}(\epsilon)$. Thus, $C_{q_i}^{\epsilon}=\emptyset$ unless
$q_i=2k$.

Next we define accordingly
\begin{align}
\omega^{\epsilon}&\triangleq \sum_i\lambda_i
\int_{C^{\epsilon}_{q_i}} \mu_i,\nonumber \\
\omega^{(\epsilon;j,a_j)}&\triangleq \sum_i\lambda_i
\int_{C^{(\epsilon;j,a_j)}_{q_i}} \mu_i.  \nonumber
\end{align}

In the limit $\epsilon\to 0^+$, we recover the whole configuration space
$C_{q_i}$, so $\omega^{(\epsilon;j,a_j)}$ becomes $\omega$.

We now associate, once and for all, an index $a_j\in\{1,2\}$ to the
$j$th crossing and construct according to
(\ref{alfai}) and (\ref{zyklon})
a $k(n-3)$-cycle in $\imbr n$
which we call
$\beta(\epsilon,\delta)$, where $\epsilon$ and $\delta$ are small
but positive.
Different choices of the parameters $\epsilon$ and $\delta$
will yield homologous cycles in $\imbr n$, as far as they do not
become too large.

Any other choice of the indices $a_j$s will also produce a homologous cycle
in force of Proposition~\ref{homseg}. We denote this new cycle by
$\beta(\epsilon,\delta) + \partial\zeta(\epsilon,\delta;j,a_j)$,
for a suitable $(k(n-3)+1)$-chain
$\zeta(\epsilon,\delta;j,a_j)$.

We now want to compute the integral
\[
I(\epsilon,\delta,j,a_j)=
\int_{\beta(\epsilon,\delta)}\omega^{(\epsilon;j,a_j)}
\]
and show that
\[
\lim_{\delta\to0^+} I(\epsilon,\delta,j,a_j) = 0.
\]
In fact, if in the explicit choice of $\beta(\epsilon,\delta)$
we have the same $a_j$ as in
$\omega^{(\epsilon;j,a_j)}$, then the above integral is zero for any $\delta$.
Otherwise, it is equal to
\[
\int_{\zeta(\epsilon,\delta;j,a_j)}
d\omega^{(\epsilon;j,a_j)}
\]
by Stokes' Theorem. The form $\omega^{(\epsilon;j,a_j)}$ is not closed.
But the main point here is that it is defined also on the space
$\immpsr 1 n (j)\subset \immsr 1 n$ given by the elements $\gamma$ for which
$\gamma(t^j_1)=\gamma(t^j_2)$ is the only double point, and $t_j^i\neq 0$.
This means that the above integral is well-defined also for
$\delta=0$ where it vanishes by dimensional reasons. But this value is
also equal to the limit $\delta\to0^+$ of the integral.

Next we consider  all
the  sets $C_{q_i}^{(\epsilon;j,a_j)}$ simultaneously. They are not disjoint,
but we can redefine
them (in an obvious way) so as to make them disjoint.
With this  we have
\[
\sum_j\int_{\beta(\epsilon,\delta)}\omega^{(\epsilon;j,a_j)}
+\int_{\beta(\epsilon,\delta)}
\omega^{\epsilon}=\int_{\beta(\epsilon,\delta)} \omega.
\]
Observe now that this expression is independent of $\delta>0$
for $\delta$ small enough; in fact,
different values of $\delta$ correspond to homologous cycles.

In particular, to compute the r.h.s.\ we can take the limit $\delta\to0^+$.
But in this limit the first term on the l.h.s.\ vanishes as proved above,
and the integral of $\omega^\epsilon$
over $\beta(\epsilon,\delta)$ is like the integral of
a Dirac-type current concentrated on the points of $C_{2k}$
whose projections on $S^1$  yield exactly
the set of those distinct
points that are in pairs identified by $\gamma\in \immpsr k n$.

Saying that only the forms on $C_{2k}$ survive is tantamount as saying
that only the chord diagrams contained in $\omega$
survive. Our final task is to prove
that, when each chord connects
two points on $S^1$ that are directly identified by $\gamma$
(i.e., that are in the pre-image of the same double point), then the
corresponding integral is $1$, otherwise it is zero.

A chord connecting two (small) intervals that contain no other vertices
can be seen as the $(n-3)$-form obtained
by integrating $v^n$ over $I\times I$ with some identifications.
To each point in $S^{n-3}$  we assign a way of lifting  one of
these two intervals
and this generates a sphere $S^{n-2}$.  The total integral associated
to this chord will then be given by the linking number of
this sphere $S^{n-2}$ with the other interval (seen as an indefinite
line). This linking number is zero if the image of the two intervals
does not contain the double point and, in force of Proposition
\ref{convicted}, is $1$ otherwise.
\end{proof}

\begin{Rem}
Observe that for $n>4$ no care has to be taken in the actual choice
of the symmetric
form $v^n$, as follows from Prop.~\ref{prop-Iexact}.
If $n=4$, then different choices of $v^n$ might yield different
homomorphisms from the graph cohomology to the de~Rham cohomology of
$\imbr 4$. However, the proof of the above proposition shows
that the evaluation of $I(v^4)(\Gamma)$,
for $\Gamma\in H^{k,0}(\mathcal{D}_e)$,
on the subspace of $H_{k(n-3)}(\imbr 4)$ spanned by the cycles
defined in \eqref{zyklon} are independent of $v^4$.
\end{Rem}

\begin{Rem}
The above proof can also be applied to the case $n=3$ assuming that one
has already made the corrections for the anomalies appearing in that case.
Notice that this proof does not reduce to the original one given in
\cite{AF}---from which it was however inspired---and seems to be a little
easier.
\end{Rem}

It is clear from the proof of Theorem~\ref{thm-nontriv} that the
image through \eqref{hzero1} of the differential form
$\omega=I(v^n)(\Lambda)$ associated to a trivalent cocycle
$\Lambda$, is equivariant in the sense of Remark~\ref{rem-eq}.
Hence, we can recast the results of the previous proof as follows:

\bth
Let $n>3$ and let
\[
\omega=I(v^n)(\Lambda)\in H^{k(n-3)}(\imbr n),
\]
with $\Lambda\in H^{k,0}(\mathcal D)$ a linear combination of trivalent
graphs.
Then the blowup map \eqref{hzero} associates to $\omega$ the linear
combinations of
the chord diagrams contained in $\Lambda$, interpreted as an element of
$H_0(\immsr k n)$. So the induced map
$H^{k,0}(\mathcal D)\to H^0(\immsr k n)$
is injective and independent of the choice of the symmetric form $v^n$.
\label{HD-HImm}
\eth

We finally remark that the proof of Theorem~\ref{thm-nontriv} also implies
the following result: a cycle of imbeddings of $S^1$ into $\bbR^n$ defined
by a special immersion is nontrivial if there is a graph cocycle in which
the chord diagram corresponding to the special immersion appears with
nonzero coefficient. We do not know however under which conditions a
chord diagram
(with no separating chords) may be part of a graph cocycle.
We will see in the next Section that, in the case of framed imbeddings
(and for $n$ odd), the situation improves and, in particular,
Theorem~\ref{tre} holds.

\subsection{Proof of Theorem~\ref{uno} and Corollary~\ref{cor-uno}}
\label{finale} The statement about the chain maps, has been
already proved in Thm.~\ref{thm-Ivn}. The other statement is a
direct consequence of the above Theorem, of
Proposition~\ref{thm-trichord} and of Corollary~\ref{crithzero}.
As for Corollary~\ref{cor-uno}, we simply recall from
subsections~\ref{ss-odd} and \ref{ss-even} that, irrespectively of
the parity, there exist cocycles of trivalent graphs in an
infinite number of orders.

\section{Framed imbeddings}
\label{sec-framed}

Let $K\colon S^1\to\bbR^n$ be an immersion. Recall that an
{\em orthonormal framing} of $K$
is a trivialization of the pulled-back bundle $K^*SO(\bbR^n)$, where
$SO(\bbR^n)\simeq\bbR^n\times SO(n)$ is the orthonormal
frame bundle of $\bbR^n$.
Equivalently we may view the framing as a map $w\colon S^1\to SO(n)$.
An orthonormal framing $w$
is said to be {\em adapted} if, for every $s\in S^1$, the
last column of $w(s)$ is equal to the normalized tangent vector
$DK(s):={\dot K(s)}/{|\dot K(s)|}$.
Let $p$ be the projection $SO(n)\to S^{n-1}$ in the quotient by $SO(n-1)$,
viewed as the subgroup of $SO(n)$ fixing the vector $(0,\dots,0,1)$.
Then saying that $w$ is a adapted framing is equivalent to stating that
the following diagram commutes:
\begin{equation}\label{frame}
\xymatrix{
 & SO(n)\ar[d]^p \\
S^1\ar[ur]^w\ar[r]_-{DK}
& S^{n-1}
}
\end{equation}
A pair $(K,w)$ consisting
of an immersion $K$ and of a adapted framing $w$ will be shortly called
a {\em framed immersion}.
We will denote by $\immrf n$ the space of framed immersions
and by $\imbrf n$ its subspace of framed imbeddings.

By forgetting of the framing we get a map $\imbrf n\to\imbr n$; so
we can pull back all cohomology classes constructed in the
previous sections. Analogously, given a framed special immersion
with $k$ transversal double points, we get a $k(n-3)$-cycle of
framed imbeddings exactly by the same construction as in
Section~\ref{vorder}. We may repeat verbatim the proof of
Theorem~\ref{thm-nontriv}, so that Theorem~\ref{HD-HImm} and
\ref{uno} and Corollary~\ref{cor-uno} immediately generalize to
the case of framed imbeddings. In the following we will show that,
for $n$ odd, there are actually more cohomology classes on $\imbrf
n$ than on $\imbr n$.

\subsection{Short chords and new cohomology classes on $\imbrf {2s+1}$}

We have seen in Section~\ref{ddco} that the boundary term
corresponding to the contraction of a short chord---viz., an
external small loop at a vertex labelled by, say, $k$---is the
form $\theta_{kk}(v^n)$, i.e., the pullback of $v^n$ by using the
composition of the projection to the $k$th point on the circle
with the normalized tangent vector map $DK$. We now wish to use a
framing, and in particular diagram \eqref{frame}, to get rid of
such terms. The plan is to find a form $\tau$ on $SO(n)$ with
$d\tau^n=p^*v^n$. \blem If $n$ is even, $[p^*v^n]$ is a nontrivial
cohomology class on $SO(n)$. If $n$ is odd, $p^*v^n$ is exact.
\elem
\begin{proof}
If $n=2s$ then the restriction map $H^*(SO(2s))\to H^*(SO(2s-1))$
is surjective and so, by Leray--Hirsch theorem, $H^*(SO(2s))$
is a free module over $H^*(S^{2s-1})$, the action of the only non trivial
generator of the cohomology of $S^{2s-1}$ on $H^*(S^{2s-1})$ being
the product by  $[p^*v^{2s}]$. Hence $p^*v^{2s}$ cannot be exact.
If $n=2s+1$, then the class of the normalized top form $v^{2s+1}$ of
$S^{2s}$ is half the Euler class and so $p^*v^{2s+1}$ is exact.
\end{proof}
Thus, we will have to restrict in the following to the case $n=2s+1$.
Using the framing, we see now that $\theta_{kk}(v^{2s+1})$ is exact;
viz., it is equal to $d\vartheta_k(\tau^{2s+1})$, where
$\vartheta_k(\tau^{2s+1})$ is
the pullback of $\tau^{2s+1}$ by the composition of the projection
to the $k$th point on the circle with the framing map $w$.

\subsubsection{Modified graph cohomology}

To keep track of the above forms, we modify the space of graphs
$\mathcal D_o^{k,m}$ by introducing a new decoration, say a cross,
that can be put on any external vertex (observe that a crossed
vertex with no edges is allowed and that there can be
more than one crossed vertex in the same graph). Crosses are
numbered and a new equivalence
relation $\Gamma\sim (-1)^\sigma\Gamma'$ is introduced, where
$\Gamma'$ is the same graph as $\Gamma$ but for the numbering of crosses
which is a permutation of the one in $\Gamma$: $\sigma$ is the order of
this permutation. Let $x$ be the number of crosses in the graph
$\Lambda$. We define the gradations
\begin{align*}
\widetilde\ord\Lambda &= e - v_i+x,  \\
\widetilde\deg\Lambda &= 2 \, e - 3 \, v_i - v_e+x.
\end{align*}
We denote by $\widetilde{\mathcal D_o^{k,m}}$ the extended
space of graphs of order $k$ and degree $m$.
We also define a map
\[
\Tilde I(v^n,\tau^{2s+1})\colon \widetilde{\mathcal D_o^{k,m}}
\to\Omega^{m+(2s-2)k}(\imbrf {2s+1}),
\]
by the same rules as in subsection~\ref{ssec-morph} plus
a new rule that associates to a cross at the vertex $k$ the form
$\vartheta_k(\tau^{2s+1})$ (the products of these odd forms being determined
by the numbering of the crosses).

We also introduce a coboundary operator $\Tilde\delta$ on
$\widetilde{\mathcal D_o^{k,m}}$; it acts on all edges, arcs and
crosses one at a time. Its action on edges and arcs is the same as
for the coboundary operator $\delta$. Its action on the cross
labelled $a$ at the vertex $i$ deletes the cross and produces an
external small loop (with half-edges ordered consistently with the
orientation of the circle); the sign associated to this operation
on a graph $\Gamma$ is defined to be
$(-1)^{\widetilde\deg\Gamma+a}$. Observe that $\Tilde\delta$
raises $\widetilde\deg$ by $1$ and leaves $\widetilde\ord$
unchanged. It is not difficult to check that $\Tilde\delta^2=0$.

It should be clear by now that the analogue of Theorem~\ref{thm-Ivn}
for framed imbeddings is the following
\begin{Prop}
\[
\Tilde I(v^n,\tau^{2s+1})\colon (\widetilde{\mathcal D_o^{k,m}}, \tilde\delta)
\to (\Omega^{m+(2s-2)k}(\imbrf {2s+1}), d)
\]
is a chain map.
\end{Prop}
\begin{proof}
The coboundary operator $\Tilde\delta$ has been defined in such a way
that $d\Tilde I(\Gamma)$ is equal to $(-1)^{\deg \tilde I(\Gamma) + 1}
\Tilde I(\Tilde\delta\Gamma)$
if one neglects hidden faces. Thus, we have to prove that hidden
faces do not contribute. But, since forms corresponding to crosses are
basic in the fibrations $\Sigma_{r,\varsigma}\to C_{q-r+1,t-\varsigma}
(\bbR^{2s+1})$,
we may rely again on the results of Appendix~\ref{hidden}.
\end{proof}

Finally, we notice that the proof of Theorem~\ref{thm-nontriv}
generalizes immediately to the present case since the forms $\vartheta$
are well-defined on immersions. So we get the following generalization
of Theorem~\ref{uno}:
\begin{Thm}
The induced homomorphism
\[
H^{k,0}_{\tilde\delta}(\widetilde{\mathcal D_o})\to H^{(2s-2)k}(\imbrf {2s+1})
\]
is injective.
\end{Thm}

\subsection{Relation between $H^{k,0}_{\tilde\delta}
(\widetilde{\mathcal D_o})$ and $H^{k,0}_{\delta}(\mathcal D_o)$}

First, we consider another graph complex, namely
our old $\mathcal D_o^{k,m}$, but with a new coboundary
operator $\underline\delta$ that differs from $\delta$ because it does
not act on arcs joining the vertices of short chords.
This new  coboundary operator  raises $\deg$ by $1$ and leaves $\ord$
unchanged.

We define a map
$\phi\colon\mathcal D_o^{k,m}\to\widetilde{\mathcal D_o^{k,m}}$.
If $\Gamma$ does not contain short chords, we set $\phi(\Gamma)=\Gamma$.
If instead $\Gamma$ contains a short chord, we proceed as follows:
First we regard $\Gamma$ as an element of $\widetilde{\mathcal D_o^{k,m}}$.
Then we pick a short chord (say, joining vertex $i$ to vertex $j$)
and replace $\Gamma$ by $\Gamma-\sigma(i,j)\Gamma'$,
where $\Gamma'$ is obtained by replacing the short chord by a crossed
vertex with label $\min\{i,j\}$, while all other vertices are renumbered
consequently. We repeat this for each short chord and number the crosses
according to this order. Clearly
 $\phi$ is well-defined (i.e., it does not depend on the order
in which we consider the short chords) and moreover we have
the following
\bprop
$\phi\colon(\mathcal D_o^{k,m},\underline\delta)
\to(\widetilde{\mathcal D_o^{k,m}},\tilde\delta)$ is a chain
map inducing an isomorphism:
\[
H^{k,0}(\phi)\colon
H^{k,0}_{\underline\delta}(\mathcal D_o)\stackrel{\sim}{\to}
H^{k,0}_{\tilde\delta}(\widetilde{\mathcal D_o}).
\]
\eprop
As a consequence we have:
\bcor
There is an injective homomorphism
\[
H^{k,0}_{\delta}(\mathcal D_o)\incul H^{k,0}_{\tilde\delta}
(\widetilde{\mathcal D_o})
\]
\ecor

\begin{proof}
In Proposition~\ref{thm-trichord} we showed that
no element of $H^{k,0}_{\delta}(\mathcal D_o)$ contains a graph with short
chords. Since $\underline\delta$ is equal to $\delta$ on elements of
$\mathcal D_o$ not containing short chords, we immediately deduce that
$H^{k,0}_{\delta}(\mathcal D_o)\subseteq
H^{k,0}_{\underline\delta}(\mathcal D_o)\cong
H^{k,0}_{\tilde\delta}(\widetilde{\mathcal D_o}).$
\end{proof}

\subsection{Relation between $H^{k,0}_{\tilde\delta}
(\widetilde{\mathcal D_o})$ and  Bar-Natan's algebra
$\mathcal A_k(\bigcirc)$}

The advantage of dealing with the complex $(\mathcal D_o^{k,m},
\underline\delta)$
instead of $(\widetilde{\mathcal D_o^{k,m}},\Tilde\delta)$ is that
$H^{k,0}_{\underline\delta}(\mathcal D_o)$ is isomorphic to the dual of the
space $\mathcal A_k(\bigcirc)$ of trivalent
graphs with a distinguished circle and oriented vertices, modulo STU
relations.

More precisely a {\em BN graph} is defined as
a connected graph made of
a distinguished oriented circle and a certain number of
edges, which are allowed to meet in two type of vertices:
internal vertices in which three edges meet, and external
vertices, in which an edge meets the distinguished circle.
External vertices are {\em  oriented} by one of the two
possible cyclic orderings of the half edges emanating from that vertex.
The degree of a BN graph is the number of its edges
minus the number of internal vertices.

$\mathcal D(\bigcirc)$ is set to be the real vector space generated by
BN graphs. $\mathcal D(\bigcirc)$ is graded by the degree
of the graphs, $\mathcal D_k(\bigcirc)$ being the component of degree $k$.

\begin{figure}[ht!]
\cl{\resizebox{8cm}{!}{\includegraphics{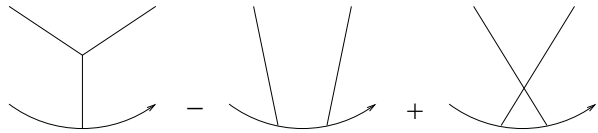}}}
\caption{$STU$ relation}
\label{fig:stu}
\end{figure}

Let $STU$ be the subspace generated by the linear combinations of
graphs as in figure \ref{fig:stu}, where the undrawn part is the same
for the three graphs that are considered.

We define $\mathcal A(\bigcirc)\triangleq\mathcal D(\bigcirc)/STU$ and
we denote its dual by $\mathcal A^*(\bigcirc)$. $\mathcal A(\bigcirc)$ is
graded and we denote the $k$th component by $\mathcal A_k(\bigcirc)$.
The space $\mathcal A^*(\bigcirc)$ is also know as the space of {\em
weight systems}.

\bprop
\label{casino}
There is an isomorphism
$\mathcal A^*_k(\bigcirc) \cong H^{k,0}_{\underline\delta}(\mathcal D_o)$
\eprop

The proof is given in Theorem 1 of ref. \cite{AF}.

\begin{Rem}
\label{sadio}
It can be shown, following
the proof of Prop.\ \ref{casino}, that given $\Gamma\in
\mathcal D_k(\bigcirc)$ and $w\in \mathcal A^*_k(\bigcirc)$ such that
$w(\Gamma)\neq 0$, then the image of $w$ in
$H^{k,0}_{\underline\delta}(\mathcal D_o)$, under the isomorphism of
Prop.\ \ref{casino}, is a cocycle  containing a graph
that, modulo decoration, is identical to $\Gamma$.
\end{Rem}

It is well known \cite{BN90, BN95}
that any finite dimensional Lie algebra with an invariant inner product
gives rise to a family of elements of $\mathcal A^*(\bigcirc)$.
In particular it is possible to
construct elements of $\mathcal A^*_k(\bigcirc)$ which are non zero on
graphs containing short chords, and hence elements of
$H^{k,0}_{\underline\delta}(\mathcal D_o)$ which contains graphs with short
chords.
This implies that $H^{k,0}_{\tilde\delta}(\widetilde{\mathcal D_o})\cong
H^{k,0}_{\underline\delta}(\mathcal D_o)$ is actually
strictly bigger than $H^{k,0}_{\delta}(\mathcal D_o)$.

\subsection{Proof of Theorem~\ref{tre}}
\label{prooftre}
For every given chord diagram there exist
cocycles of zero degree in $(\mathcal D_o^{k,m},\Hat\delta)$
containing it. This follows from Remark \ref{sadio} and from the fact that
the weight system determined
by the Lie algebra $\mathfrak{gl}(n)$ is nonzero on every chord
diagram.

\section{On the cohomology of the spaces of imbeddings and immersions of
$S^1$ into $\reali^n$} \label{chen}

In this Section we will relate the cohomology of the spaces
$\imbr{n}$ and $\immr{n}$ of imbeddings and immersions of $S^1$
into $\reali^n$. This will allow us to consider examples of
nontrivial cohomology classes of $\imbr{n}$ not coming from
trivalent graphs and even not coming from our graphs. Moreover
these last classes are not equivariant w.r.t. the action of
$\mathrm{Diff}^+(S^1)$.

We will first review known facts about the computation of
$H^*(\immr{n})$ by Chen's iterated integrals, and then we will
consider the restriction of the above cohomology classes to
$\imbr{n}$.


Given a Riemannian simply connected manifold $M$, we naturally
associate to any $\psi\in \imm{M}$ its normalized derivative
$D\psi$, which is an element of the loop space $\mathcal{L}SM$,
where $SM$ denotes the tangent sphere bundle. We denote by
$\mathcal{L} X$ the  free loop space over a manifold $X$.

Smale proved \cite{Smale} (see also \cite{Bry}) the following: \bth
The map $D:\imm{M}\to \mathcal{L}SM$ is a weak homotopy
equivalence. \eth In particular the space $\immr n$ is weakly
homotopy equivalent to $\mathcal{L}S^{n-1}$.

The cohomology of the loop space of a simply connected manifold
$X$ can be computed via Chen's iterated integrals \cite{Chen} (for
a purely algebraic approach see, for example,\cite{NdoTho}).

Let us consider the $m$-dimensional simplex $\Sigma_m:=\{
(t_1,\ldots t_m) \in \reali^m : \,\,\, 0 \leq t_1 \leq \ldots \leq
t_m \leq 1 \}$ and the map $\Phi$
\[
\Phi\colon
\begin{array}[t]{ccc}
\Sigma_m \times \mathcal{L}X & \lora &X\times \dots \times X \\
((t_1,\ldots ,t_m),\gamma) & \rsa & (\gamma(t_1),\dots,
\gamma(t_m)).
\end{array}
\]
{}From a collection of forms $w_1,\ldots ,w_m$ in $\Omega^*(X)$,
we obtain the form $pr^\ast_1\,w_1 \wedge,$ $\ldots, \wedge
pr^\ast_m\,w_m$ on $X^m$, where $pr_i:X^m\to X$ denotes the
projection onto the $i$th factor. The iterated Chen integral is
defined by pullback via the map $\Phi$ and subsequent integration
on $\Sigma_m$. The result is a form on $\mathcal{L}X$,

Moreover, any form $w_0\in\Omega^*(X)$ defines also a form of the
same degree on $\mathcal{L}X$, denoted again as $w_0$,  obtained
by pullback via the evaluation map at the initial point of the
loop. In summary, we consider  the following  class of forms on
$\mathcal{L}X$ with degree $\sum_{i=0}^{m} \deg(w_i) -m$:
\begin{equation}
w_0P(w_1,\cdots,w_m) \triangleq w_0\int_{\Sigma_m}
\Phi^\ast (pr^\ast_1\,w_1 \wedge,\ldots \wedge
pr^\ast_m\,w_m). \label{phast}
\end{equation}
When all the forms $w_i$ are equal to a given form $w$, we use
the shorter notation $w_0P_m(w)$.

When $M=\R^n$, it is sufficient to consider $w_i=v^n$, where $v^n$
is a fixed normalized top form on $S^{n-1}$. In fact, we have
\cite{GJP}:

\bprop
The algebra $H^*\left(\mathcal{L}S^{n-1},\reali\right)$ is
the associative and graded commutative algebra with the following
set of generators and relations:
\begin{itemize}
\item
For $n$ odd, the generators are $P_1(v^n)P_{2s}(v^n)$ in degree
$(2s+1)(n-2))$ and $v^nP_{2s}(v^n)$ in degree $(2s+1)(n-2)+1$,
with $s=0,1,\dots$. The relations are given by the requirement
that the product of any two generators is zero.
\item
For $n$ even, the generators are $v^n$ in degree $n-1$ and
$P_1(v^n)$ in degree $n-2$ with the only relation $(v^n)^2=0$.
\end{itemize}
\eprop

This determines uniquely the real cohomology of $\immr n$, which
is obtained by pulling back to $\immr n$ the above classes via the
map $D$ (normalized derivative).

Some of the above cohomology classes are equivariant
with respect to the action of the group $\diff^+(S^1)$ of
orientation preserving diffeomorphisms of $S^1$, on $\immr n$.
More precisely, one can prove \cite{GJP} that irrespectively of
the dimension $n,$ the forms $v^nP_k(v^n)$ do not yield equivariant
classes on $\immr n$ and that the equivariant generators are given by
$D^*P_1(v^n)$ for $n$ even and by $D^*P_1(v^n)D^*P_{2s}(v_n)$,
with relations as before, for $n$ odd.

\begin{Rem}\label{rem-equiv-forms}
The equivariant forms $D^*P_k(v^n)$ can be constructed using the morphisms
$I(v^n)$ of equations \eqref{oddmorph} and \eqref{evenmorph}. Let us
consider in fact a graph $\Gamma_k\in \mathcal D^{k,k}$ consisting
of $k$ small loops attached to $k$ different external vertices.
Since $\Gamma_k$ contains no chord, $I(v^n)(\Gamma_k)$ is a well
defined differential form on $\immr n$, and we have an (obvious)
identification $D^*P_k(v^n) = I(v^n)(\Gamma_k)$. For instance, the
closure of $D^*P_k(v^n)$ can be seen as a consequence of the equation
$\delta\Gamma_k=0$, which holds both in the odd and in the even case.
\end{Rem}

We now consider the restriction of the above
forms and classes on $\immr n$ to forms and classes on $\imbr n$.
It is convenient to assume that $v^n$ is a symmetric form
(see Def~\ref{def-gf}). Then we define $\theta(v^n)$ to be the
configuration space integral $I(v^n)(\Gamma'_1)$ where $\Gamma'_1
\in\mathcal D^{1,0}$ is the graph with two external vertices
joined by one chord, and $I(v^n)$ the homomorphisms of
eq.~\eqref{oddmorph} or ~\eqref{evenmorph}, depending on the parity of
$n$. Clearly $\theta(v^n)$ is zero if $n$ is even because the
corresponding graph is zero by the
relation \eqref{relazp}. On the contrary, if $n$ is odd, one has the
following:
\blem
If $n$ is odd, we have
\[
d\theta(v^n)=-2D^*P_1(v^n).
\]
\label{thetap}
\elem

\begin{proof}
First, let us write
\[
\theta(v^n) = \int_{S^1\times I} \phi^*v^n
\]
where
\begin{equation}
\phi\colon
\begin{array}[t]{ccc}
S^1\times I \times \imbr{n} &\to& S^{n-1} \label{fi} \\
(s,t,\psi) &\rsa& \frac{\psi(s+t)-\psi(s)}{||\psi(s+t)-\psi(s)||,}
\end{array}
\end{equation}
We have
\begin{equation}
d\theta(v^n)=\int_{S^1\times \partial I}i^* \phi^*v^n= D^*P_1(\alpha^*v^n)-
D^*P_1(v^n),
\label{dtheta}
\end{equation}
where the map $i$ denotes the restriction to the
boundary $S^1\times\partial I$ and $\alpha$ is the antipodal map on
$S^{n-1}$.
\end{proof}

We now prove the following:
\bth
If $n$ is odd, then the inclusion of $\imbr n$ in
$\immr n$ induces the zero map in (real) cohomology.
\label{nodd}
\eth

\proof
The restriction of $D^*P_1(v^n)$ to $\imbr n$ is exact in force of
Lemma \ref{thetap}. When we restrict the generators
$D^*[P_1(v^n)P_{2s}(v^n)]$ in degree $(2s+1)(n-2))$ to $\imbr n$,
we obtain trivial cocycles, as shown in the following equation:
\begin{multline}
-2D^*[P_1(v^n)P_{2s}(v^n)]=\\
=d\left(\theta(v^n) D^*P_{2s}(v^n)- (1/2)\theta^2(v^n)
D^*[v^nP_{2s-2}(v^n)]\right).\label{Pk-exact}
\end{multline}

Next we  prove that also the generators $D^*[v^nP_{2s}(v^n)]$
represent trivial cocycles when restricted to  $\imbr n$.

First we notice that $D^*v^n$ is exact when restricted to $\imbr
n$. Consider in fact the map $\widehat\phi:I\times \imbr n\to
S^{n-1}$ given by the restriction of the map $\phi$ of equation
(\ref{fi}) to $\{0\}\times I\times \imbr n\subset S^1\times
I\times \imbr n$ and define
\begin{equation}
  \label{defTheta}
  \Theta(v^n)\triangleq \int_I\widehat\phi^* v^n.
\end{equation}
We conclude that
\[
d\Theta(v^n) = -2D^* v^n.
\]
Finally, we have
$$
-2D^*[v^nP_{2s}(v^n)]=d\left(\Theta(v^n)D^*P_{2s}(v^n)-
(1/2)\Theta^2(v^n) D^*P_{2s-1}(v^n)\right).
\eqno{\qed}$$

Thanks to the identification $D^*P_k(v^n) = I(v^n)(\Gamma_k)$ of
Remark~\ref{rem-equiv-forms}, Lemma~\ref{thetap} and
equation~\eqref{Pk-exact} can be seen in graph cohomology
as the existence of $\Gamma'_k\in
\mathcal D^{k,k-1}_o$ such that $\delta_o(\Gamma'_k)=\Gamma_k$.

Theorem~\ref{nodd} has been proved in \cite{BT}, but only for the
case $H^1(\immr 3)$.

If $n$ is even, then we know that $\theta(v^n)$ is zero and
moreover it is easy to see that $\Theta(v^n)$ (as defined in
\eqref{defTheta}) is closed. Thus, if $n$ is even, we cannot
conclude anymore that the inclusion $\imbr n\incul \immr n$
induces the zero morphism in cohomology. On the contrary, we find
that the classes $D^*P_1(v^n)$, $D^*v^n$ and $D^*[v^n P_1(v^n)]$
remain nontrivial when restricted to $\imbr n$. We first have:

\bprop If $n$ is even, then the restriction of $D^*[v^n P_1(v^n)]$
to\break $\imbr n$ is a nontrivial cocycle of degree $2n-3$.
\label{p1veven} \eprop

\proof
Let us consider the real Stiefel manifold
\[
V_{n,2}=SO(n)/SO(n-2).
\]
It is an orientable compact manifold of dimension $(2n-3)$. To
each $(\mathbf x, \mathbf y)\in V_{n,2}$ we associate the
following element of  $\imbr n$:
\begin{equation}
\label{cycleP1v} t\rsa ((1-\cos t)\;\mathbf x  + \sin t\;\mathbf
y).
\end{equation}
The integral of the cocycle $D^*[v^n P_1(v^n)]$ over the image of
the map \eqref{cycleP1v} is given by the degree of the map:
\[
\begin{array}[t]{ccc}
f:S^1 \times V_{n,2}(\reali)\to S^{n-1}\times S^{n-1}\\
(t, \mathbf{x}, \mathbf{y})\rsa \left(\mathbf{y}, \sin
t\;\mathbf{x}+ \cos t\;\mathbf{y}\right).
\end{array}
\]
It is evident that
\begin{multline}
\nonumber f^{-1}
((0,\cdots,0,0,1),(0,\cdots,0,1,0))=\\
=(\pi/2,(0,\cdots,0,1,0),(0,\cdots,0,0,1))\,\cup\\
\cup\,(3\pi/2,(0,\cdots,0,-1,0),(0,\cdots,0,0,1))
\end{multline}
and
$$
\deg f=1+(-1)^{n}.
\eqno{\qed}$$

\bcor If $n$ is even, then the restrictions of $D^*P_1(v^n)$ and
$D^*v^n$ to  $\imbr n$ are nontrivial cocycles of degree $n-2$
and, respectively, $n-1$. \label{p1-vn-even} \ecor
\begin{proof}
$D^*[v^n P_1(v^n)]$ is the product of $D^*v^n$ and $D^*P_1(v^n)$.
\end{proof}

For completeness, we explicitly construct examples of cycles
of $\imbr n$ below, nontrivially paired with the above cocycles.

If $n$ is even, then a nontrivial element $\Gamma_{n-2}\in
H_{n-2}(\imbr n)$ is given by the map that assigns to each
normalized vector $\mathbf{x}\in S^{n-2}$, represented as
$(x_1,x_2,\cdots, x_{n-1},0)\in \reali^{n}$, the imbedded loop in
$\reali^n$ given by
\[
t\rsa (x_1 (1-\cos t), x_2 (1-\cos t), \dots,x_{n-1} (1-\cos t),
\sin t).
\]
By computing the degree of the normalized derivative of this map,
it turns out that the evaluation of $D^*P_1(v^n)$ on
$\Gamma_{n-2}$ is $-2$.

Again, for $n=2s$, a nontrivial element $\Gamma_{n-1}\in
H_{n-1}(\imbr n)$ is given by the map that assigns to each
normalized vector
\[
\mathbf x= (x_1,y_1,x_2,y_2,\cdots,x_s,y_s) \in S^{2s-1}\subset
\reali^{2s}
\]
the loop based at $\mathbf x$ given by
\begin{multline}  \nonumber
t\rsa
 (x_1 (1-\cos t) + y_1 \sin t, -x_1 \sin t + y_1 (1-\cos t),
\cdots\\
\cdots, x_s (1-\cos t) + y_s \sin t, -x_s \sin t + y_s (1-\cos
t)),
\end{multline}
and the evaluation of $D^*v^n$ on $\Gamma_{n-1}$ yields
$(-1)^{s}$.

\begin{Rem}
Corollary~\ref{p1-vn-even} extends Thm.~\ref{uno} since it
provides an example in which the morphism~\eqref{evenmorph} is
nontrivial on a nontrivalent graph cocycle $\Gamma_1\in H^{1,1}
(\mathcal D_e)$.
\end{Rem}

\begin{Rem}
To complete the discussion on the nature of the restriction to
$\imbr n$ of nontrivial classes on $\immr n$, we should also say
something about $D^*P_k(v^n)$ and $D^*[v^n P_k(v^n)]$ when $n$ is
even and $k>1$.

At the moment, however, we have no cycles to pair with them in
order to prove that they are nonzero nor have an argument to prove
their vanishing. We may only comment that these classes are images
under \eqref{evenmorph} of nontrivial cocycles of graphs (in the
sense of Section~\ref{ddco}). Therefore, triviality of the above
classes would imply that \eqref{evenmorph} is not injective.
\end{Rem}

\begin{Rem}
We have chosen $v^n$ to be a given symmetric form. Notice however
that the main results of this Section (Theorem~\ref{nodd},
Proposition~\ref{p1veven} and Corollary~\ref{p1-vn-even}) are independent
of the choice of the volume form on $S^{n-1}$).
\end{Rem}

\appendix\small
\section{Appendix: Integrals along the boundary and vanishing theorems}
\label{hidden}

The goal of this Appendix is to show that, as long as $n>3$ and we
use only symmetric forms, {\em hidden faces}
do not contribute to Stokes' theorem.

\subsection{Codimension-one faces}
Let us begin with a description of
the codimension-one faces of $\de C_{q,t}(\reali^n)$
These faces can be divided into three classes:
\begin{description}
\item[Type I] $s\ge2$ of the points in $\reali^n$
collapse together;
\item[Type II] $s\ge1$ of the points in $\reali^n$ escape
together to infinity;
\item[Type III] $r\ge1$ of the points in $S^1$ and $s\ge0$
of the points in $\reali^n$, with $r+s\ge2$, collapse together; we
denote these faces by the symbol $\Sigma_{r,s}$.
\end{description}

Boundary faces of {\bf type I} or {\bf type II}
are given
by $\widehat{C}_s(\reali^n)\times C_{q,t-s+1}(\reali^n)$
and $\widehat{C}_{s+1}(\reali^n)$ $\times C_{q,t-s}(\reali^n)$ respectively.
Here $\widehat{C}_k(\reali^n)$ is obtained from $(\reali^n)^k/G$ ---
$G$ being the group of global translations and scalings --- by
blowing up all diagonals.
Each $\widehat{C}_k(\reali^n)$ is a compact manifold
with corners. In the simplest case, $k=2$, we have $\widehat{C}_2(\reali^n)
=S^{n-1}$. Notice however that in the following we will only need to consider
the interior of $\widehat{C}_k(\reali^n)$ which is easily identified with
$C^0_{k-1}(\reali^n)/\reali_+$, where the group $\reali_+$ acts by rescaling
all components.

A boundary face $\Sigma_{r,0}$ is a fibration over
$C_{q-r+1,t}(\reali^n)$ got by pulling back $\widehat{C}_r(TS^1)$.
In the case when $s>0$, the description of $\Sigma_{r,s}$ is a little bit
longer. Actually, we will only describe the interior of $\Sigma_{r,s}$
(which is enough for Stokes' theorem). Let us first define
\[
\begin{split}
S^{n-1}\ltimes \mathcal{C}_{r,s} \triangleq
\{
a\in S^{n-1},\ (x_1,\dots,x_r)\in \reali^r,\
(y_1,\dots,y_s)\in(\reali^n)^s \mid \\
a(x_i)\not=y_j,\ \forall i,j;\
x_i=x_j \Rightarrow i=j;\
y_i=y_j \Rightarrow i=j
\},
\end{split}
\]
where $a\in S^{n-1}$ is seen as a linear map $\reali\to\reali^n$.

Next introduce the equivalence relations
\begin{align*}
(a;x_1,\dots;y_1,\dotsc) &\sim
(a;\lambda\,x_1,\dots;\lambda\,y_1,\dotsc),\\
(a;x_1,\dots;y_1,\dotsc) &\sim
(a;x_1+\xi,\dots;y_1+a(\xi),\dotsc),
\end{align*}
with $\lambda>0$ and $\xi\in\reali$, and
denote by $S^{n-1}\ltimes\widehat{C}_{r,s}$
the quotient. This is a fibration over $S^{n-1}$ whose fiber
will be denoted by $\widehat{C}_{r,s}$.

Then the interior $\mbox{int}(\Sigma_{r,s})$
of a {\bf type III} face, corresponding to the
collapse of $r$
points on $S^1$ together with $s$  points in $\reali^n$,
is a fibration over $C_{q-r+1,t-s}(\reali^n)$
obtained by pulling back
$(S^{n-1}\ltimes\widehat{C}_{r,s})$
via the composition $\widehat D$ of the normalized derivative
$D:\imbr n\to \mathcal
LS^{n-1}$ with the evaluation at the position on
$S^1$ where the points have collapsed:
\[
\begin{array}{ccc}
\mbox{int}(\Sigma_{r,s})& \lora & S^{n-1}\ltimes\widehat{C}_{r,s}\\
\downarrow & & \downarrow \nonumber \\
C_{q-r+1,t-s} (\reali^n) & \stackrel{\widehat D}{\lora} & S^{n-1}
\nonumber
\end{array}
\]

\begin{Rem}
What is really important for all the subsequent considerations are
the dimensions of the fibers of the above spaces. A simple computation
yields
\begin{equation}
\begin{split}
d_{0,s} &\dot= \dim\widehat{C}_s(\reali^n) = ns - n - 1,\\
d_{r,s} &\dot= \dim\widehat{C}_{r,s} = r + ns - 2  \quad\text{for } r>0.
\end{split}\label{drs}
\end{equation}
\end{Rem}

We end with the following
\bdefiniz
A codimension-one face is called a {\em principal face}\/ if it
is of type I with $s=2$, or of type II with $s=1$, or
of type III with $r+s=2$.
All other codimension-one faces are called {\em hidden faces}.
\edefiniz

\subsection{Subgraphs corresponding to codimension-one faces}
In Section~\ref{ddco}, we have explained how to associate integrals
of tautological forms over configuration spaces to decorated graphs.
Here we want to represent the integration
along codimension-one faces in terms of subgraphs.

First we introduce the notion of {\em decorated graphs with a
point at infinity}. These are just decorated graphs with one extra
vertex labelled by $\infty$ and the condition that $\infty$ is not
the end-point of any edge. Of course, there is a one-to-one
correspondence $\Gamma\mapsto\Tilde\Gamma$ between decorated
graphs and decorated graphs with a point at infinity.

A subgraph $\Gamma'$ of a graph $\Tilde\Gamma$ is a graph
whose vertices are a subset of the vertices in $\Tilde\Gamma$ and
whose edges and arcs are the subset of the edges and arcs
of $\Tilde\Gamma$
whose end-points are distinct and both in the subset of vertices
of $\Gamma'$.
(Notice in particular that an external small loop based at a vertex
in $\Gamma'$ is not considered as an edge of $\Gamma'$.)

A subgraph $\Gamma'$ of a decorated graph with a point at infinity
is said
to be {\em admissible} if one of the following conditions is fulfilled:
\begin{description}
\item[Type I] $\Gamma'$ has $s\ge2$ internal vertices and no other vertices;
\item[Type II] $\Gamma'$ has $s\ge1$ internal vertices plus the
$\infty$-vertex;
\item[Type III] $\Gamma'$ has $r\ge1$ external vertices and
$s\ge0$ internal vertices, with $r+s\ge2$.
\end{description}
The correspondence between admissible subgraphs and
codimension-one faces should be clear.

\bdefiniz
An admissible subgraph will be called principal (hidden) if
the corresponding face is principal (hidden).
\edefiniz

If $\Gamma'$ is an admissible subgraph of $\Gamma$, we define the
reduction of $\Gamma$ modulo $\Gamma'$, which we will denote by
$[\Gamma/\Gamma']$, as follows. First we define two vertices in
$\Gamma$ to be equivalent if they both belong to $\Gamma'$. Then
the set of vertices of $[\Gamma/\Gamma']$ is the set of vertices
of $\Gamma$ modulo the equivalence relations, while the set of
edges and arcs contains all the edges and arcs of $\Gamma$ that do
not belong to $\Gamma'$. The vertex in $[\Gamma/\Gamma']$ whose
pre-image in $\Gamma$ is $\Gamma'$ will be denoted by
$v(\Gamma',\Gamma)$.

Given any graph $\Gamma$, we will denote by $\alpha_\Gamma$ the
corresponding form (product of tautological forms associated to the edges).

Following the description of the preceding Section, we can compute
the form on the space of imbeddings obtained by integrating
$\alpha_\Gamma$ along the codimension-one face corresponding to
$\Gamma'$ as follows:  we integrate along the configuration space
corresponding to $[\Gamma/\Gamma']$  the product of
$\alpha_{[\Gamma/\Gamma']}$ and a form which we denote by
$\alpha(\Gamma/\Gamma')$.

For faces of type I, $\alpha(\Gamma/\Gamma')$ is the form
$I(\Gamma')$ obtained by integrating $\alpha_{\Gamma'}$ along
$\Gamma'$, and as such vanishes unless it is a zero-form.

For faces of type II, sending a subgraph $\Gamma'$ 
to infinity is, roughly speaking, like shrinking to a point the 
complementary subgraph. More precisely, we consider the subgraph
$\widetilde\Gamma$ of $\Gamma$ consisting of the vertices which
are not in $\Gamma'$, and of those arcs and edges whose end-points
do not belong to $\Gamma'$. Then we form the quotient
$\Gamma''=\Gamma/\widetilde\Gamma$ and denote by $v$ the vertex of
$\Gamma''$ whose pre-image is $\widetilde\Gamma$. Now
$\alpha(\Gamma/\Gamma')$ is obtained by integrating
$\alpha_{\Gamma''}$ over the subspace of the configuration space
corresponding to $\Gamma''$, in which $v$ is fixed to $0\in\R^n$.
Again, we have that $\alpha(\Gamma/\Gamma')$ vanishes unless it is
a zero-form.

For a face of type III, $\alpha(\Gamma/\Gamma')$ is the pullback
of $I(\Gamma')$ via the map $D$ (normalized derivative), combined
with the evaluation at the point corresponding to the vertex
$v(\Gamma',\Gamma)$. As a consequence, $I(\Gamma')$ and
$\alpha(\Gamma/\Gamma')$ vanish unless their degree is less or
equal to $n-1$.

\subsection{Vanishing theorems}
Now we can state the following:
\bth
If\/ $\Gamma'$ is a principal subgraph, then $I(\Gamma')$
vanishes unless
\begin{enumerate}
\item $\Gamma'$ is of type I  and contains exactly one
edge,
or
\item $\Gamma'$ is of type III with $r=s=1$ and contains exactly one
edge,
or
\item $\Gamma'$ is of type III with $r=2$ and has at most one edge.
\end{enumerate}
In the first two cases or in the third case with no  edges
$I(\Gamma')=\pm1$. In the third case with one edge,
$I(\Gamma')=\pm v^n$. \label{thm-princ}
\eth

\begin{proof}
Since $\Gamma'$ is principal, it contains at most two vertices.
Thus, by our definition of graphs, there is at most
one edge in $\Gamma'$ if it of type I or III, and no edge if $\Gamma'$
is of type II.

If there are no edges, then there is nothing to be
integrated. However, the fiber corresponding to $\Gamma'$ has dimension
strictly positive (equal to $n-1$) if it is of type I or II, or if it is
of type III with $r=s=1$.

If however the face is of type III with $r=2$
the fiber dimension is zero, and we get $I(\Gamma')=1$.

If there is exactly one edge and $\Gamma'$ is of type I or it is of type III
with $r=s=1$, we have to
integrate the $(n-1)$-dimensional form $v$ over the fiber
which is $S^{n-1}$.
This yields $I(\Gamma')=\pm1$.

If $\Gamma'$ is of type III with $r=2$ and contains exactly one   edge
(to which corresponds $v^n$), then the fiber is zero-dimensional,
and $I(\Gamma')=\pm v^n$.
\end{proof}

\begin{Rem} The cases in which $I(\Gamma')$ does not vanish correspond
exactly to the integrations over principal faces considered in
Section~\ref{ddco}. The sign $I(\Gamma')=\pm 1$ is the sign $\sigma$
of  (\ref{sigma1}), (\ref{sigma2}) and (\ref{sigmaedge}), i.e.,
the sign associated to the operator $\delta$; the case
$I(\Gamma')= \pm v^n$ corresponds to the contraction of an arc
joining the end-points of a short chord. \label{rem-princ}
\end{Rem}

As for the hidden faces, we have the following result, whose proof
is contained in the next subsubsection.
\bth
If\/ $\Gamma'$ is an
admissible subgraph corresponding to a hidden face and $n>3$, then
$\alpha(\Gamma/\Gamma')=0$.
\label{thm-vt}
\eth

\subsubsection{Proof of Theorem~\ref{thm-vt}}
We follow the strategy proposed in \cite{K} and \cite{Thu}.

Throughout this subsection, by valence of a vertex we mean the number of
edges ending at that vertex.

We start by proving the following basic facts.
\blem
If\/ $\Gamma'$ is hidden and contains a zero-valent vertex different
from $\infty$, then $\alpha(\Gamma/\Gamma')=0$.
\label{lem-1}
\elem
\blem
If\/ $\Gamma'$ is hidden and contains a univalent internal vertex,
then $\alpha(\Gamma/\Gamma')=0$.
\label{lem-2}
\elem
\blem
If\/ $\Gamma'$ is hidden and contains a bivalent internal vertex,
then $\alpha(\Gamma/\Gamma')=0$.
\label{lem-3}
\elem
\begin{proof}[Proof of Lemma~\ref{lem-1}]
Let $i$ be the $0$-valent vertex and $x_i$ the corresponding point.
We perform the $x_i$-integration first. Let us compute the dimension
of this fiber. Since $\Gamma'$ is hidden, we can use
the other points to fix translations and scalings (or just scalings
if $\Gamma'$ is of type II). So the fiber is
$n$-dimensional if $i$ is internal and $1$-dimensional if $i$ is external.

The form to be integrated over this fiber is however of degree zero
since no   edge ends at $i$.
\end{proof}
\begin{proof}[Proof of Lemma~\ref{lem-2}]
As in the previous proof, we perform the $x_i$-integration
first, where $i$ is an internal univalent vertex now. As before,
the fiber is $n$-dimensional.

The form to be integrated over this fiber is the tautological
form corresponding to the edge ending at $i$. Thus, it has degree $n-1$.
\end{proof}
\begin{proof}[Proof of Lemma~\ref{lem-3}]
Now we perform the $n$-dimensional integration over $x_i$, where
$i$ is an internal bivalent vertex.

Let $j$ and $k$ be the end-points of the two   edges starting at $i$.
If $j=k$, then the integrand form is $\theta_{ij}^2=0$.

So assume $j\not=k$. The $x_i$-integration can be the extended to
$\reali^n$ with only $x_j$ and $x_k$ blown up. Let us denote by
$C(\reali^n;x_j,x_k)$ this space.
We must then compute
\[
\alpha = \int_{C(\reali^n;x_j,x_k)} \theta_{ij}\,\theta_{ik}.
\]
Following Kontsevich \cite{K}, we consider
the involution $\chi$ of $C(\reali^n;x_j,x_k)$ that
maps $x_i$ to $x_j+x_k-x_i$. Observe that $\chi$ is orientation
preserving (reversing) if $n$ is even (odd). Moreover,
\[
\begin{split}
\chi^*\theta_{ij} &= \theta_{ki} =(-1)^n\,\theta_{ik},\\
\chi^*\theta_{ik} &= \theta_{ji} =(-1)^n\,\theta_{ij}.\
\end{split}
\]
Thus,
\[
\alpha = (-1)^n \int_{C(\reali^n;x_j,x_k)} \theta_{ik}\,\theta_{ij}=
(-1)^n\,(-1)^{n-1}\int_{C(\reali^n;x_j,x_k)} \theta_{ij}\,\theta_{ik}=
-\alpha.
\]
Therefore, $\alpha=0$.
\end{proof}

Thanks to the above Lemmas, we can from now on assume that all internal
vertices in $\Gamma'$ are at least trivalent and that all external
vertices are at least univalent.
As a consequence, we have
\[
2e' \ge r+3s,
\]
where $e'$ is the number of edges in $\Gamma'$, $r$ is the number of external
vertices (possibly zero) and $s$ that of internal vertices.
Therefore, since any edge carries an $(n-1)$-form,
\[
\deg \alpha(\Gamma/\Gamma') \ge
\begin{cases}
(n-1)\frac{r+3s}2 - d_{r,s} & \text{if $\Gamma'$ is of type I or III},\\
(n-1)\frac{3(s+1)}2 - d_{0,s+1} & \text{if $\Gamma'$ is of type II},
\end{cases}
\]
with $d_{r,s}$ defined in \eqref{drs}.
In order to prove that $\alpha(\Gamma/\Gamma')$ vanishes, it is then enough
to check that $\deg \alpha(\Gamma/\Gamma')>0$ for $\Gamma'$ of type I or II,
and that $\deg \alpha(\Gamma/\Gamma')>n-1$ for $\Gamma'$ of type III (otherwise
the form to be integrated over the fiber vanishes by dimensional reasons).
Finally, a straightforward computation yields
\begin{equation}
\deg \alpha(\Gamma/\Gamma') \ge
\begin{cases}
\frac{(n-3)s}2 + n + 1 & \text{if $\Gamma'$ is of type I},\\
\frac{(n- 3)s}2 + \frac{5n-1}2 & \text{if $\Gamma'$ is of type II},\\
\frac{(n-3)(r+s-2)}2 + n - 1 & \text{if $\Gamma'$ is of type III}.
\end{cases}   \label{degIG}
\end{equation}
The result now follows from the facts that $\Gamma'$ is hidden (so $r+s>2$)
and that $n>3$.
\begin{Rem}
In the case $n=3$, the previous argument fails for $\Gamma'$ of type III,
and suitable corrections have to be found (see \cite{BT, Thu, Poi}).

In the case $n=2$, the argument fails in general and more refined techniques
have to be considered (see \cite{K3}).
\end{Rem}

\subsection{Forms depending on a parameter}
As in the proof of Prop.~\ref{prop-Iexact}, we now
consider tautological forms that are pullbacks
of a form $\Tilde v^n$ in $\Omega^{n-1}(S^{n-1}\times[0,1])$
which is closed and symmetric and integrates to one on the sphere for any value
of the parameter $t\in[0,1]$. The corresponding configuration space integral
$\Tilde I(\Gamma)$ yields then a form on $\imbr n\times[0,1]$.
It is not difficult to check that Thm.~\ref{thm-princ} still holds with
$I$ replaced by $\Tilde I$. Moreover, we have the following
generalization of Thm.~\ref{thm-vt}:
\bth
If\/ $\Gamma'$ is an admissible subgraph corresponding to a hidden face
and
$n>4$, then $\Tilde I(\Gamma')=0$.
\label{vt01}
\eth
\begin{proof}
We proceed as in the proof of Thm.~\ref{thm-vt} till the paragraph just before
\eqref{degIG}.

Now we observe that, by dimensional reasons, $\Tilde I(\Gamma')$ vanishes if
$\deg \Tilde I(\Gamma')>1$ for $\Gamma'$ of type I or II,
or $\deg\Tilde I(\Gamma')>n$ for $\Gamma'$ of type III.

But these inequalities are again satisfied by \eqref{degIG}.
\end{proof}
\begin{Rem}
In the case $n=4$, the above result might not hold.
For example, consider the 5-dimensional face $\Sigma_{3,1}$ and assume
that in $\Gamma'$ there are three edges (so, a 9-form). After integration,
we are left with a 4-form on $S^3\times[0,1]$, and there is no apparent
reason why it should vanish.

In particular, consider the cocycle of figure~\ref{figura3}.
When all the vertices in the second graph collapse, we are exactly in the
same situation as above. Since there are no other edges, we can now perform
the integration over $[0,1]$. We are finally left with the integral along
the knot of the pullback via the normalized derivative of a top form
$w^4$ on $S^3$.
Thus, on the r.h.s.\ of \eqref{Iv10n}, we must add a multiple (possibly zero)
of $P_1(w^4)$, that is, of the lowest cohomology class obtained by
restriction from the space of immersions.
\label{v4}
\end{Rem}

\normalsize\vglue -35pt

\hbox{}

\Addresses\recd

\end{document}